\documentclass[preprint]{elsarticle} 

\makeatletter
\def\ps@pprintTitle{%
\let\@oddhead\@empty
\let\@evenhead\@empty
\def\@oddfoot{\centerline{\thepage}}%
\let\@evenfoot\@oddfoot}
\makeatother

\usepackage{etoolbox}
\newcommand{\preprintdate}{July 16, 2026}
\patchcmd{\MaketitleBox}{\footnotesize\itshape\elsaddress\par\vskip36pt}{\footnotesize\itshape\elsaddress\par\parbox[b][36pt]{\linewidth}{\vfill\hfill\textnormal{\preprintdate}\hfill\null\vfill}}{}{}%
\patchcmd{\pprintMaketitle}{\footnotesize\itshape\elsaddress\par\vskip36pt}{\footnotesize\itshape\elsaddress\par\parbox[b][36pt]{\linewidth}{\vfill\hfill\textnormal{\preprintdate}\hfill\null\vfill}}{}{}%

\usepackage[
colorlinks=true,%
breaklinks=true,%
linkcolor=blue,
urlcolor=blue,%
citecolor=blue,%
pdftitle={Gradient-enhanced spline dimensional decomposition
  for uncertainty quantification with limited  training samples}, 
pdfkeywords={Spline dimensional decomposition, 
  Uncertainty quantification, Gradient-enhanced surrogate modeling, 
  Energy-balanced scaling, Ridge regression, Orthonormal B-splines},	
pdfauthor={Eunho Heo; Dongjin Lee},		
bookmarksopen=false,
pdfpagemode=UseNone]{hyperref}
\usepackage{multicol}
\usepackage{algorithm}
\usepackage{algorithmicx}
\usepackage{algpseudocode}
\usepackage{tabularx}
\usepackage{array}
\usepackage{multirow}
\usepackage{threeparttable}
\usepackage{booktabs}
\usepackage{graphicx}
\usepackage{arydshln}
\usepackage{setspace}   
\usepackage{xcolor}
\usepackage[mathlines,running]{lineno}
\usepackage{epstopdf}
\usepackage{graphicx}
\usepackage[numbers]{natbib}
\usepackage{algorithm}
\usepackage{algpseudocode}
\usepackage{booktabs}
\usepackage{placeins}
\usepackage{float}
\usepackage{multirow}

\newcommand{\Astack}{\mathbf{A}_{\mathrm{stack}}}
\newcommand{\bstack}{\mathbf{b}_{\mathrm{stack}}}
\newcommand{\Wrow}{\mathbf{W}_{\mathrm{row}}}
\newcommand{\Atil}{\tilde{\mathbf{A}}_{\mathrm{stack}}}
\newcommand{\btil}{\tilde{\mathbf{b}}_{\mathrm{stack}}}
\newcommand{\Dpen}{\mathbf{D}}

\usepackage[margin = 1.2in]{geometry}
\usepackage[T1]{fontenc}
\usepackage[english]{babel}
\usepackage[latin1]{inputenc}
\usepackage{mathtools}
\usepackage{amsmath}
\usepackage{amsfonts}
\usepackage{amsthm}
\usepackage{amssymb}
\usepackage{array}
\usepackage{subcaption}
\usepackage{siunitx}

\usepackage{color}
\usepackage{url}
\usepackage{cleveref}
\usepackage{graphicx}
\usepackage{breakcites}
\usepackage{soul} 
\usepackage{enumitem}
\usepackage[title,titletoc,toc]{appendix}

\newcommand{\RNum}[1]{\uppercase\expandafter{\romannumeral #1\relax}}
\usepackage{amsthm}
\makeatletter
\def\els@aparagraph[#1]#2{\elsparagraph[#1]{#2\@addpunct{.}}}
\def\els@bparagraph#1{\elsparagraph*{#1\@addpunct{.}}}
\makeatother
\newtheorem{remark}{Remark}
\newtheorem{assumption}{Assumption}

\newtheorem{problem}{Problem}

{\noindent {\textbf{Proof}:} }%
{\hfill $\Box$ \\[1ex] }

\newcommand{\bit}{\begin{itemize}}
\newcommand{\eit}{\end{itemize}}
\newcommand{\ben}{\begin{enumerate}}
\newcommand{\een}{\end{enumerate}}

\graphicspath{{./figures/}}
\begin{document}	
\begin{frontmatter}

\title{Gradient-enhanced spline dimensional decomposition
  for uncertainty quantification with limited  training samples}
		
		
\author[HYU1]{Eunho Heo}
\ead{gjrkd426@hanyang.ac.kr}
\author[HYU2]{Dongjin Lee\corref{cor1}}
\ead{dlee46@hanyang.ac.kr}

\cortext[cor1]{Corresponding author}
\address[HYU1]{Department of Automotive Engineering (Automotive-Computer Convergence), \\ 
Hanyang University, Seoul, South Korea}
\address[HYU2]{Department of Automotive Engineering, Hanyang University, Seoul, South Korea}

\begin{abstract}
A spline dimensional decomposition (SDD) surrogate effectively
represents high-dimensional engineering responses with localized features
and complex nonlinearities in uncertainty quantification (UQ).
However, limited training data can make coefficient estimation from function values severely ill-conditioned.
We propose gradient-enhanced SDD (GE-SDD), which trains the surrogate using function values and partial derivatives. A diagonal row-weight matrix balances the function and derivative blocks by their Frobenius norms. We solve the balanced system through ridge regression in probability-weighted Sobolev coordinates and select the regularization parameter using grouped $K$-fold cross-validation to prevent information leakage. Mapping the solution back to the $L^2$-orthonormal SDD basis preserves closed-form mean and variance estimates.

We evaluate the proposed GE-SDD on a two-dimensional continuous exponential function, a linear dynamical system with three uncertain parameters, and a 30-dimensional 25-bar truss. GE-SDD is more accurate than standard SDD and uses gradients more robustly than gradient-enhanced Kriging. GE-SDD achieves a median NRMSE of $1.022\%$ on the nonsmooth benchmark, compared with $8.731\%$ for Kriging. For the truss, GE-SDD yields lower NRMSE and more accurate standard-deviation estimates than Kriging at moderate training sizes and above. Overall, the benefits of gradient augmentation depend on input dimension, basis resolution, training size, and the target UQ quantity.
\end{abstract}


\end{frontmatter}


\section{Introduction}
\label{sec:introduction}

Uncertainty in material properties, geometry, and loads 
propagates to engineering predictions, making uncertainty quantification (UQ) essential for reliable analysis and design~\cite{
Xiu2002,noh2011reliability}.
Forward UQ estimates output statistics, such as the
mean and variance, by propagating input uncertainty through a computational model.
Monte Carlo simulation (MCS) is broadly applicable and 
asymptotically unbiased, but its
slow convergence becomes prohibitively expensive 
when each sample requires a
finite-element or computational fluid dynamics
simulation.

Surrogate models reduce this cost by approximating the simulator from 
a limited set of training samples. Common surrogates include
spectral expansions, such as polynomial chaos expansion (PCE)~\cite{Xiu2002,novak2021variance} and polynomial dimensional decomposition (PDD)~\cite{Rahman2008,Rahman2011}; collocation and 
sparse-grid methods~\cite{Babuska2007,Gerstner1998,Xiu2005}; 
Gaussian process regression (GPR)~\cite{Bouhlel2019}; reduced-order models~\cite{mcbane2021component} and neural networks~\cite{ryu2024physics}. Multifidelity strategies further reduce the sampling cost 
by adaptively combining surrogates and models of varying 
fidelity~\cite{chaudhuri2018multifidelity,lee2023multifidelity}.
Neural networks offer considerable flexibility in approximating
highly nonlinear responses through richly parameterized
function classes defined by their weights and biases.
However, they often require many training samples and 
iterative optimization. They also lack the closed-form 
moment estimates available from measure-consistent orthogonal expansions. 
PCE and PDD compute the mean and variance directly from their 
coefficients~\cite{Rahman2008,Rahman2011}, but their global polynomial bases often struggle 
with localized nonlinearities, steep gradients, and near-discontinuities common in structural dynamics, fracture, and multiphysics.

Spline dimensional decomposition (SDD) addresses this
limitation by combining an
ANOVA decomposition with locally supported B-spline basis
functions~\cite{Rahman2022,audoux2020nurbs}. Knot vectors partition each input domain, allowing SDD to capture local variations that are difficult to resolve with PCE or PDD~\cite{Rahman2021,
Jahanbin2022,Dixler2021}.
The hierarchical ANOVA structure and orthonormal basis also provide analytical mean and variance estimates without additional sampling. 
Combined with score functions, SDD yields 
closed-form second-moment sensitivities for robust design optimization~\cite{Lee2022} and supports stochastic isogeometric analysis on multipatch
domains~\cite{Jahanbin2022}.
However, standard SDD estimates its
 coefficients from function values alone.
In strongly nonlinear problems with high-dimensional inputs, the number of basis functions can exceed the available samples and make coefficient estimation ill-conditioned.

Gradient observations can improve 
coefficient estimation when training data are limited.
Intrusive sensitivity methods differentiate the governing equations with
respect to input parameters and solve the resulting sensitivity or adjoint equations. Applications include structural shape
sensitivity and
adjoint-based aerodynamic
design~\cite{Giles2000}.
These methods are effective when user can access and differentiate the governing equations, but they require internal simulator code or dedicated sensitivity solvers. Such requirements limit their use with complex or legacy black-box models. 

Non-intrusive surrogates can use gradients obtained from analytical sensitivities, adjoint
solvers, automatic differentiation, finite differences,
or built-in sensitivity modules.
Gradient-enhanced Kriging
(GE-Kriging), formulated within a co-Kriging
framework, jointly models function values and partial derivatives through auto- and cross-covariances~\cite{Morris1993,
Bouhlel2019}.
For a sufficiently smooth kernel, 
differentiating the kernel analytically provides the required covariance terms. 
GE-Kriging has produced more accurate predictions than conventional Kriging for fixed sampling budgets and has reduced the number of samples needed to reach a target accuracy in problems with up to 100 input dimensions and engineering problems
with up to 15 dimensions~\cite{Bouhlel2019}. However, its computational cost grows
cubically with the total number of function and derivative observations, while maximum likelihood hyperparameter estimation becomes increasingly
expensive as the input dimension increases. 

Gradient-enhanced polynomial methods instead convert function values and derivatives into linear equations for coefficient estimation. At each training sample, one function value and one partial derivative per input dimension provides $1+d$ equations for a $d$-dimensional problem~\cite{Peng2016}. Gradient-augmented compressed sensing stacks these equations and applies sparsity-promoting regularization, such as the Least Absolute Shrinkage and Selection Operator (LASSO)~\cite{Candes2006,
Tibshirani1996}. Derivative observations can therefore reduce the number of training samples required for reliable sparse polynomial approximation, particularly in UQ problems with high-dimensional inputs.

Regularized regression has also improved dimensional-decomposition surrogates without gradient data. A sparsity-promoting
D-MORPH regression method first computes a sparse LASSO solution and then embeds that solution within D-MORPH regression to train PDD surrogates for global sensitivity analysis~\cite{Lee2024}.
This approach improves coefficient estimation with limited data but does not use derivative observations.

Despite these developments, no
gradient-enhanced SDD method has used derivative observations within a Sobolev-type
coefficient-estimation framework, to the authors' knowledge. Such a formulation must address two practical challenges. 
First, the function and derivative blocks of the design matrix can differ in Frobenius norm by several orders of magnitude across input dimensions. Directly stacking these blocks allows high-norm blocks to dominate the regression loss. Second, each function value and its derivatives come from the same training sample. Cross-validation must therefore keep these observations in the same fold; otherwise, gradient information can leak into the validation set and produce overly optimistic estimates. 

To overcome these difficulties, we propose a
gradient-enhanced SDD (GE-SDD) method for UQ with limited
training samples.
The main contributions are as follows:
\begin{enumerate}
\item[(1)] We formulate GE-SDD as a probability-weighted Sobolev-type regression that combines function values and partial derivatives. Each training sample contributes additional deterministic equations while preserving the standard SDD representation for moment estimation.
\item[(2)] We develop a block-wise, energy-balanced scaling strategy that prevents function and derivative blocks with large Frobenius norms from dominating the regression.
\item[(3)] We map the coefficients estimated in Sobolev coordinates back to the original $L^2$-orthonormal SDD coordinates. This mapping preserves closed-form mean and variance estimates without additional Monte Carlo sampling.
\item[(4)] We use grouped $K$-fold cross-validation that keeps all observations from each training sample in the same fold, preventing information leakage during regularization-parameter selection.
\end{enumerate}
The remainder of this paper is organized as follows.
Section~\ref{sec:background} introduces the UQ problem, SDD framework, and Sobolev-space formulation. Section~\ref{sec:gesdd} presents GE-SDD. Section~\ref{sec:numerical} evaluates the method on three benchmark problems of increasing complexity. Section~\ref{sec:conclusions} summarizes the main findings.

\section{Background and problem formulation}
\label{sec:background}

This section provides the mathematical foundation
for GE-SDD. 
Section~\ref{subsec:notation} defines the notation, and Sections~\ref{subsec:variables} and~\ref{subsec:uq-objective} formulate the random variables and UQ objective. Section~\ref{subsec:sdd} introduces the ANOVA structure, orthonormal B-spline basis, and truncated SDD approximation. Finally, Section~\ref{subsec:sobolev} develops the Sobolev-space setting for the gradient-enhanced formulation in Section~\ref{sec:gesdd}.

\subsection{Mathematical notation}
\label{subsec:notation}

Let $\mathbb{N}$, $\mathbb{N}_0$, $\mathbb{R}$, and
$\mathbb{R}_0^+$ denote the sets of positive integers,
non-negative integers, real numbers, and non-negative
real numbers, respectively.
For $N\in\mathbb{N}$, define the bounded rectangular input domain 
$\mathbb{A}^N
  :=\prod_{k=1}^N \mathbb{A}_k,
  ~
  \mathbb{A}_{k}:=[a_k,b_k]\subset\mathbb{R},
  ~ k=1,\ldots,N.
  \label{eq:input-domain}$
For a nonempty subset
$u\subseteq\{1,\ldots,N\}$ with cardinality $|u|$,
denote by $\mathbf{x}_u$ and $\mathbf{X}_u$ the
subvectors of $\mathbf{x}$ and $\mathbf{X}$ indexed by
$u$, respectively, and define
$\mathbb{A}^u:=\prod_{k\in u}\mathbb{A}_{k}$.
We denote expectation and variance by
$\mathbb{E}[\cdot]$ and $\mathbb{V}\mathrm{ar}[\cdot]$,
respectively, with respect to the input 
probability distribution defined below. 

\subsection{Input and output random variables}
\label{subsec:variables}
Let $(\Omega,\mathcal{F},\mathbb{P})$ be an abstract
probability space, where $\Omega$ is the sample space,
$\mathcal{F}$ is a $\sigma$-algebra on $\Omega$, and
$\mathbb{P}:\mathcal{F}\to[0,1]$ is a probability
measure.
Consider an $N$-dimensional random vector
$\mathbf{X}:=(X_1,\ldots,X_N)^\top:\Omega\to \mathbb{A}^N,$
representing the uncertain inputs of the computational
model. The joint cumulative distribution function (CDF) of
$\mathbf{X}$ is $F_{\mathbf{X}}(\mathbf{x}):=\mathbb{P}\!\left[\bigcap_{k=1}^{N}\{X_k\le x_k\} \right].$
When $F_{\mathbf{X}}$ is absolutely continuous, the
corresponding joint probability density function (PDF)
is
$f_{\mathbf{X}}(\mathbf{x}):=\partial^N F_{\mathbf{X}}(\mathbf{x})/(\partial x_1\cdots\partial x_N)$.
The probability measure induced by $\mathbf{X}$ on its
input domain is represented by
$(\mathbb{A}^N,\mathcal{B}^N,f_{\mathbf{X}}(\mathbf{x})\,\mathrm{d}\mathbf{x})$,
where $\mathcal{B}^N:=\mathcal{B}(\mathbb{A}^N)$ is the Borel
$\sigma$-algebra on $\mathbb{A}^N$.

\begin{assumption}
\label{assum:1}
The random vector $\mathbf{X}$ satisfies the following
conditions.
\begin{enumerate}
  \item[(1)] The components $X_k$, $k=1,\ldots,N$, are
        mutually independent but not necessarily
        identically distributed.
  \item[(2)] Each $X_k$ is defined on a bounded interval
        $\mathcal{A}_{(k)}=[a_k,b_k]\subset\mathbb{R}$.
  \item[(3)] Each $X_k$ has an absolutely continuous
        marginal distribution function $F_{X_k}$ and a
        marginal probability density function $f_{X_k}$
        supported on $[a_k,b_k]$.
\end{enumerate}
\end{assumption}
Under Assumption~\ref{assum:1}, the joint density
factorizes as $f_{\mathbf{X}}(\mathbf{x})=
\prod_{k=1}^{N}f_{X_k}(x_k).$
The bounded support also ensures the existence of the
finite moments required for constructing the
measure-consistent orthonormal spline basis.

Let $y:\mathbb{A}^N\to\mathbb{R}$
denote a deterministic computational model or
quantity-of-interest (QoI) mapping.
The corresponding output random variable is defined by
the composition $Y:=y(\mathbf{X}):\Omega\to\mathbb{R}$.
The computational model is assumed to be square
integrable with respect to the input probability measure, such that
$
  y\in L^2(\mathbb{A}^N,\mathcal{B}^N,f_{\mathbf{X}}(\mathbf{x})\mathrm{d}\mathbf{x}),
 ~
  \mathbb{E}[Y^2]
  =
  \int_{\mathbb{A}^N}
  y(\mathbf{x})^2
  f_{\mathbf{X}}(\mathbf{x})\,
  \mathrm{d}\mathbf{x}
  <\infty.
  \label{eq:L2}
$
%
\subsection{Objectives of uncertainty quantification}
\label{subsec:uq-objective}
We draw the training samples from the input probability measure or select them using a deterministic experimental design. 
The goal is to construct a surrogate that accurately approximates $y(\mathbf{X})$ over $\mathbb{A}^N$ and estimates the mean and variance of 
$Y=y(\mathbf{X})$, such that 
\begin{align}
  \mu_Y
  :=
  \mathbb{E}[Y]
  =
  \int_{\mathbb{A}^N}
  y(\mathbf{x})
  f_{\mathbf{X}}(\mathbf{x})\,
  \mathrm{d}\mathbf{x}\qquad\text{and}\qquad 
   \sigma_Y^2
  :=
  \mathbb{V}\mathrm{ar}[Y]
  =
  \int_{\mathbb{A}^N}
  \bigl(y(\mathbf{x})-\mu_Y\bigr)^2
  f_{\mathbf{X}}(\mathbf{x})\,
  \mathrm{d}\mathbf{x}
\end{align}
from limited function-value observations.
In engineering applications, evaluating
$y(\mathbf{x})$ may require an expensive finite-element,
computational fluid dynamics, or multiphysics simulation.
The number of available training samples is
therefore often limited.
Let $M$ denote the number of training samples, each providing 
one function-value observation: 
$
  \left\{
    \left(
      \mathbf{x}^{(i)},
      y(\mathbf{x}^{(i)})
    \right)
  \right\}_{i=1}^{M},
  ~
  \mathbf{x}^{(i)}\in \mathbb{A}^N.
  \label{eq:function-training-data}
$


In the gradient-enhanced setting, 
each of the $M$ training samples also provides
partial derivatives with respect to the input
variables.
The resulting training set is 
\begin{equation}
  \left\{ \begin{bmatrix} \mathbf{x}^{(i)} \\ y(\mathbf{x}^{(i)}) \\ \nabla y(\mathbf{x}^{(i)}) \end{bmatrix} \right\}_{i=1}^M \subset \mathbb{R}^{2N+1}.
  \label{eq:ge-training-data}
\end{equation}
We treat the derivatives as local constraints associated with the same $M$ training points, rather than as independent samples. Analytical sensitivities, adjoint methods, automatic differentiation, finite differences, or dedicated sensitivity modules can provide the required derivatives.

\subsection{Spline dimensional decomposition}
\label{subsec:sdd}

SDD uses the analysis-of-variance (ANOVA)
dimensional decomposition of
$y(\mathbf{X})$~\cite{Hoeffding1948,Rabitz1999}.
Under Assumption~\ref{assum:1}, the response has the
unique representation
$
  y(\mathbf{X})
  =
  y_\emptyset
  +
  \sum_{\emptyset\ne u\subseteq\{1,\ldots,N\}}
  y_u(\mathbf{X}_u),
$
where $y_\emptyset:=\mathbb{E}[y(\mathbf{X})]$ and each
component function $y_u$ depends only on $\mathbf{X}_u$.
For each $k\in u$, $y_u$ satisfies
$
  \int_{a_k}^{b_k}
  y_u(\mathbf{x}_u)
  f_{X_k}(x_k)\,
  \mathrm{d}x_k
  =
  0
$
for almost every $\mathbf{x}_{u\setminus\{k\}}$.
This zero-mean condition ensures the uniqueness and mutual
orthogonality of the ANOVA components.

SDD approximates each nonconstant component using locally
supported, measure-consistent orthonormal B-spline basis
functions. For each input coordinate $X_k$, let $p_k$,
$\boldsymbol{\xi}_k$, and $n_k$ denote the spline degree,
open knot vector, and number of univariate basis functions,
respectively. The corresponding orthonormal basis functions
are denoted by
$\psi^k_{i_k,p_k,\boldsymbol{\xi}_k}$,
$i_k=1,\ldots,n_k$, where
$\psi^k_{1,p_k,\boldsymbol{\xi}_k}\equiv 1$ and the
remaining functions have zero mean under the probability
measure of $X_k$. For a non-empty subset $u$, the associated
tensor-product basis function is
\begin{equation}
  \Psi^u_{\mathbf{i}_u,\mathbf{p}_u,\boldsymbol{\Xi}_u}
  (\mathbf{x}_u)
  :=
  \prod_{k\in u}
  \psi^k_{i_k,p_k,\boldsymbol{\xi}_k}(x_k).
  \label{eq:tp-basis}
\end{equation}
The detailed construction and orthonormalization procedure
is provided in Appendix~\ref{app:basis-construction}.

For each non-empty interaction subset $u$, define the
reduced multi-index set
$\bar{\mathcal{I}}_{u,\mathbf{n}_u}
  :=
  \{
    \mathbf{i}_u= (i_{k_1}, \dots, i_{k_{|u|}}) :
    2\le i_k\le n_k,\; k\in u
  \}.
  \label{eq:reduced-index}$
Consequently, every retained tensor-product basis
function has zero mean and contributes only to the
non-constant part of the ANOVA decomposition.
For
$\mathbf{i}_u\in\bar{\mathcal{I}}_{u,\mathbf{n}_u}$,
the retained tensor-product basis functions satisfy
\begin{equation}
  \mathbb{E}\!\left[
    \Psi^u_{\mathbf{i}_u,\mathbf{p}_u,\boldsymbol{\Xi}_u}
    (\mathbf{X}_u)
  \right]
  =0 .
  \label{eq:retained-zero-mean}
\end{equation}
Under Assumption~\ref{assum:1}, for
$\mathbf{i}_u\in\bar{\mathcal{I}}_{u,\mathbf{n}_u}$ and
$\mathbf{j}_v\in\bar{\mathcal{I}}_{v,\mathbf{n}_v}$, the
retained tensor-product basis functions satisfy
\begin{equation}
  \mathbb{E}\!\left[
    \Psi^u_{\mathbf{i}_u,\mathbf{p}_u,\boldsymbol{\Xi}_u}
    (\mathbf{X}_u)
    \Psi^v_{\mathbf{j}_v,\mathbf{p}_v,\boldsymbol{\Xi}_v}
    (\mathbf{X}_v)
  \right]
  =
  \begin{cases}
    1,
      & u=v
        \text{ and }
        \mathbf{i}_u=\mathbf{j}_v,
        \\[2mm]
    0,
      & \text{otherwise}.
  \end{cases}
  \label{eq:ortho-multi}
\end{equation}

For a maximum interaction order $S\le N$, the total
number of retained SDD basis functions, including the
constant term, is
\begin{equation}
  L_{S,\mathbf{p},\boldsymbol{\Xi}}
  =
  1+
  \sum_{s=1}^{S}
  \sum_{\substack{
        u\subseteq\{1,\ldots,N\}\\
        |u|=s}}
  \prod_{k\in u}(n_k-1).
  \label{eq:basis-count}
\end{equation}

\paragraph{Truncated SDD approximation.}
Using the retained tensor-product basis functions, the
theoretical $S$-variate SDD representation is
\begin{equation}
  y_{S,\mathbf{p},\boldsymbol{\Xi}}(\mathbf{X})
  :=
  c_\emptyset
  +
  \sum_{\substack{
        \emptyset\ne u\subseteq\{1,\ldots,N\}\\
        1\le |u|\le S}}
  \sum_{\mathbf{i}_u\in
        \bar{\mathcal{I}}_{u,\mathbf{n}_u}}
  c^u_{\mathbf{i}_u,\mathbf{p}_u,\boldsymbol{\Xi}_u}
  \Psi^u_{\mathbf{i}_u,\mathbf{p}_u,\boldsymbol{\Xi}_u}
  (\mathbf{X}_u),
  \label{eq:sdd}
\end{equation}
where $c_\emptyset:=y_\emptyset=\mathbb{E}[Y]$
is the exact constant coefficient, and
$c^u_{\mathbf{i}_u,\mathbf{p}_u,\boldsymbol{\Xi}_u}
\in\mathbb{R}$ are the non-constant SDD expansion
coefficients. In data-driven SDD construction, all
coefficients are estimated from the available observations.
The estimated constant coefficient, denoted by
$\hat c_\emptyset$, therefore approximates the exact ANOVA
constant $y_\emptyset$.

Owing to the orthonormality and zero-mean properties of
the retained basis functions, the mean and variance of
the theoretical truncated SDD representation are

\begin{equation}
  \mathbb{E}\!\left[
    y_{S,\mathbf{p},\boldsymbol{\Xi}}(\mathbf{X})
  \right]
  =
  c_\emptyset
  =
  y_\emptyset,
  \qquad
  \operatorname{Var}\!\left[
    y_{S,\mathbf{p},\boldsymbol{\Xi}}(\mathbf{X})
  \right]
  =
  \sum_{\substack{
        \emptyset\ne u\subseteq\{1,\ldots,N\}\\
        1\le |u|\le S}}
  \sum_{\mathbf{i}_u\in
        \bar{\mathcal{I}}_{u,\mathbf{n}_u}}
  \left(
    c^u_{\mathbf{i}_u,\mathbf{p}_u,\boldsymbol{\Xi}_u}
  \right)^2.
  \label{eq:sdd-moment}
\end{equation}
The same formulas apply to the estimated coefficients, so the fitted SDD provides the mean and variance without additional Monte Carlo sampling.

\subsection{Sobolev framework for GE-SDD}
\label{subsec:sobolev}

GE-SDD augments standard SDD
regression with first-order partial derivatives
of the retained basis functions.
These basis functions must therefore have
square-integrable weak first derivatives.
The response $y$ must also provide the corresponding
partial derivatives at the training samples with gradient observations.
Throughout this subsection, let $w:\mathbb{A}^N\to\mathbb{R}$ be Lebesgue-measurable. 
The space of 
square-integrable functions is $L^2(\mathbb{A}^N)
  :=
  \left\{
    w:\mathbb{A}^N\to\mathbb{R}
    \;\text{measurable}:
    \int_{\mathbb{A}^N} w(\mathbf{x})^2\,\mathrm{d}\mathbf{x}
    <\infty
  \right\}.
  \label{eq:L2-lebesgue}$
For any $w_1,w_2\in L^2(\mathbb{A}^N)$, the corresponding inner product is
$\langle w_1,w_2\rangle_{L^2(\mathbb{A}^N)}
:=\int_{\mathbb{A}^N} w_1(\mathbf{x})w_2(\mathbf{x})\,
\mathrm{d}\mathbf{x}$.
This Lebesgue-measure space differs from the probability-weighted space
$L^2(\mathbb{A}^N,\mathcal{B}^{N},f_{\mathbf{X}}(\mathbf{x})\mathrm{d}\mathbf{x}$) defined
in~\eqref{eq:L2}.

We denote by $C_c^\infty(\mathbb{A}^N)$ the space of infinitely
differentiable test functions
$\varphi:\mathbb{A}^N\to\mathbb{R}$ whose compact support lies in the interior of $\mathbb{A}^N$.
For each $k\in\{1,\ldots,N\}$, the weak partial derivative
of $w$ with respect to $x_k$ is a function
$\partial w/\partial x_k\in L^2(\mathbb{A}^N)$ satisfying, if
\begin{equation}
  \int_{\mathbb{A}^N}
  w(\mathbf{x})
  \frac{\partial\varphi}{\partial x_k}(\mathbf{x})\,
  \mathrm{d}\mathbf{x}
  =
  -
  \int_{\mathbb{A}^N}
  \frac{\partial w}{\partial x_k}(\mathbf{x})
  \varphi(\mathbf{x})\,
  \mathrm{d}\mathbf{x}
  \quad
  \text{for every }
  \varphi\in C_c^\infty(\mathbb{A}^N),
  \label{eq:weak-grad}
\end{equation}
The weak gradient is
$\nabla w:=(\partial w/\partial x_1,\ldots,
\partial w/\partial x_N)^\top$.

The first-order Sobolev space on $\mathbb{A}^N$ is then defined
explicitly as
\begin{equation}
  H^1(\mathbb{A}^N)
  :=
  \left\{
    w\in L^2(\mathbb{A}^N):
    \frac{\partial w}{\partial x_k}\in L^2(\mathbb{A}^N),\;
    k=1,\ldots,N
  \right\}.
  \label{eq:H1-def}
\end{equation}
where $\partial w/\partial x_k$ denotes the weak partial
derivative defined in~\eqref{eq:weak-grad}.
For $w_1,w_2\in H^1(\mathbb{A}^N)$, 
the standard inner product is
\begin{equation}
  \langle w_1,w_2\rangle_{H^1(\mathbb{A}^N)}
  :=
  \langle w_1,w_2\rangle_{L^2(\mathbb{A}^N)}
  +
  \sum_{k=1}^{N}
  \left\langle
    \frac{\partial w_1}{\partial x_k},
    \frac{\partial w_2}{\partial x_k}
  \right\rangle_{L^2(\mathbb{A}^N)}.
  \label{eq:H1-inner}
\end{equation}
To consider the input probability measure, we define
the probability-weighted Sobolev-type bilinear form
\begin{align}
  a_f(w_1,w_2)
  &:=
  \mathbb{E}\!\left[
    w_1(\mathbf{X})w_2(\mathbf{X})
  \right] +
  \sum_{k=1}^{N}
  \mathbb{E}\!\left[
    \frac{\partial w_1}{\partial x_k}(\mathbf{X})
    \frac{\partial w_2}{\partial x_k}(\mathbf{X})
  \right].
  \label{eq:weighted-H1-bilinear}
\end{align}
All expectations in~\eqref{eq:weighted-H1-bilinear} are finite
for bounded functions with bounded weak first derivatives; in
particular, they are finite for the piecewise-polynomial
B-spline basis functions used in this work.

B-splines of degree $p_k\ge1$ with
open knot vectors and simple interior knots are continuous piecewise polynomials with square-integrable weak first derivatives. Thus, each univariate orthonormal basis functions~\eqref{eq:ortho-univariate} and the retained
tensor-product basis functions belong to $H^1(\mathbb{A}^N)$.
For $p_k\ge2$, the univariate basis functions are $C_c^{p_k-1}$-continuous across simple interior knots. 

The Sobolev framework therefore supports the inclusion of basis derivatives in coefficient estimation.
Section~\ref{subsec:stacked} uses the probability-weighted derivative terms to construct the coordinate transformation, while Section~\ref{subsec:scaling} addresses numerical imbalances between the function-value and derivative blocks.

\begin{problem}
\label{prob:main}
Given $M$ training samples with function-value and partial-derivative observations, construct an
accurate GE-SDD surrogate in the limited-sample regime ($M < L_{S,\mathbf{p},\boldsymbol{\Xi}}$) while 
preserving coefficient-based estimates of the output 
mean and variance.
\end{problem}

\section{Gradient-enhanced spline dimensional
decomposition}
\label{sec:gesdd}

\begin{figure}
  \centering
  \includegraphics[width=.9\textwidth]
  {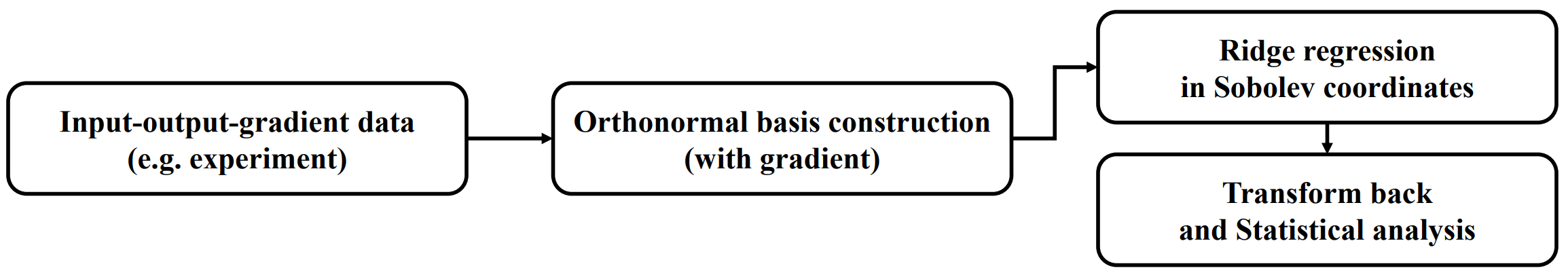}
  \caption{Flowchart for the proposed GE-SDD framework. The method transforms the $L^2$-orthonormal SDD basis into Sobolev coordinates and estimates the expansion coefficients from function and gradient observations using ridge regression.}
  \label{fig:flowchart}
\end{figure}

This section develops the proposed GE-SDD method.
Section~\ref{subsec:stacked} formulates the
probability-weighted Sobolev-transformed regression system
using both function-value and partial-derivative
observations.
Section~\ref{subsec:scaling} introduces block-wise,
energy-balanced scaling to balance the Frobenius
norms of the function and derivative blocks.
Section~\ref{subsec:ridge-cv} presents the ridge
regression and grouped $K$-fold
cross-validation for selecting the
regularization parameter.
Section~\ref{subsec:moments} presents coefficient estimation in the original $L^2$-orthonormal SDD coordinates, retaining closed-form mean and variance estimates. 
Figure~\ref{fig:flowchart} summarizes the overall workflow.

\subsection{Sobolev-transformed stacked regression system}

\label{subsec:stacked}
Denote the number of retained SDD basis functions by
$
  n_\psi
  :=
  L_{S,\mathbf{p},\boldsymbol{\Xi}}
  \label{eq:npsi-definition}
$. 
Collecting the constant and retained
tensor-product basis functions constructed in
Section~\ref{subsec:sdd} gives the multivariate
$L^2$-orthonormal SDD basis vector as
\begin{equation}
  \boldsymbol{\psi}_{S,\mathbf{p},\boldsymbol{\Xi}}
  (\mathbf{X})
  :=
  \bigl(
    \psi_1(\mathbf{X}),
    \psi_2(\mathbf{X}),
    \ldots,
    \psi_{n_\psi}(\mathbf{X})
  \bigr)^\top,
  \qquad
  \psi_1(\mathbf{X})\equiv1.
  \label{eq:basis-full}
\end{equation}
We omit the subscripts
$(S,\mathbf{p},\boldsymbol{\Xi})$ when
no ambiguity arises. Orthonormality gives
\begin{equation}
  \mathbb{E}\!\left[
    \boldsymbol{\psi}(\mathbf{X})
  \right]
  =
  (1,0,\ldots,0)^\top,
  \qquad
  \mathbb{E}\!\left[
    \boldsymbol{\psi}(\mathbf{X})
    \boldsymbol{\psi}^{\top}(\mathbf{X})
  \right]
  =
  \mathbf{I}_{n_\psi}.
  \label{eq:psi-l2-orthonormal}
\end{equation}

We obtain the first-order derivatives of the basis functions analytically using the B-spline differentiation
rule.
The extended basis matrix collects the function basis and its partial derivatives:
\begin{equation}
  \mathbf{A}(\mathbf{X})
  :=
  \begin{bmatrix}
    \boldsymbol{\psi}(\mathbf{X})
    &
    \dfrac{\partial\boldsymbol{\psi}}{\partial X_1}(\mathbf{X})
    &
    \cdots
    &
    \dfrac{\partial\boldsymbol{\psi}}{\partial X_N}(\mathbf{X})
  \end{bmatrix}
  \in
  \mathbb{R}^{n_\psi\times(1+N)}.
  \label{eq:A-extended}
\end{equation}
The first column of $\mathbf{A}(\mathbf{X})$ contains the
original SDD basis functions, and remaining
columns contain their partial derivatives with respect to the input variables.

Using the bilinear form
introduced in~\eqref{eq:weighted-H1-bilinear}, define
the derivative-augmented spline moment matrix  by
\begin{align}
  \mathbf{G}_{\mathrm{GE}}
  &:=
  \mathbb{E}\!\left[
    \mathbf{A}(\mathbf{X})
    \mathbf{A}^{\top}(\mathbf{X})
  \right]
  \nonumber\\
  &=
  \mathbb{E}\!\left[
    \boldsymbol{\psi}(\mathbf{X})
    \boldsymbol{\psi}^{\top}(\mathbf{X})
  \right]
  +
  \sum_{k=1}^{N}
  \mathbb{E}\!\left[
    \dfrac{\partial\boldsymbol{\psi}}{\partial X_k}(\mathbf{X})
    \dfrac{\partial\boldsymbol{\psi}^{\top}}{\partial X_k}(\mathbf{X})
  \right]
  \nonumber\\
  &=
  \mathbf{I}_{n_\psi}
  +
  \sum_{k=1}^{N}
  \mathbb{E}\!\left[
    \dfrac{\partial\boldsymbol{\psi}}{\partial X_k}(\mathbf{X})
    \dfrac{\partial\boldsymbol{\psi}^{\top}}{\partial X_k}(\mathbf{X})
  \right]
  \in
  \mathbb{R}^{n_\psi\times n_\psi}.
  \label{eq:gram-ge}
\end{align}

The derivative contribution in~\eqref{eq:gram-ge}
is positive semidefinite, and the function-value
contribution is $\mathbf{I}_{n_\psi}$. Therefore,
$\mathbf{G}_{\mathrm{GE}}$ is symmetric positive definite and admits the Cholesky factorization
$  \mathbf{G}_{\mathrm{GE}}
  =
  \mathbf{Q}\mathbf{Q}^{\top},
  \label{eq:gram-ge-cholesky}$
where $\mathbf{Q}$ is nonsingular and lower-triangular.
The mathematical formulation requires no diagonal regularization.
Numerical implementations may add a machine-scale diagonal perturbation to mitigate finite-precision
roundoff; this safeguard is not part of the proposed transformation.

\begin{remark}
[Assembly of the Sobolev spline moment matrix]
\label{rem:gram-assembly}
We assemble $\mathbf{G}_{\mathrm{GE}}$ over the retained
ANOVA interaction subsets.
For retained basis functions
$\Psi^u_{\mathbf{i}_u}$ and
$\Psi^v_{\mathbf{j}_v}$, the derivative term associated 
with $X_k$ can be nonzero only when
$k\in u\cap v$ and the remaining active factors satisfy
the corresponding $L^2$-orthogonality conditions. We therefore assemble $\mathbf{G}_{\mathrm{GE}}$ entrywise or
blockwise over the retained SDD index set, avoiding a full $N$-dimensional tensor-product construction.
\end{remark}

The Sobolev-transformed extended basis matrix is 
\begin{equation}
  \mathbf{B}(\mathbf{X})
  :=
  \mathbf{Q}^{-1}\mathbf{A}(\mathbf{X})
  =
  \begin{bmatrix}
    \boldsymbol{\Phi}_1(\mathbf{X})
    &
    \boldsymbol{\Phi}_2(\mathbf{X})
    &
    \cdots
    &
    \boldsymbol{\Phi}_{1+N}(\mathbf{X})
  \end{bmatrix},
  \label{eq:ge-ortho-basis}
\end{equation}
where
\begin{equation}
  \boldsymbol{\Phi}_1(\mathbf{X})
  =
  \mathbf{Q}^{-1}\boldsymbol{\psi}(\mathbf{X}),
  \qquad
  \boldsymbol{\Phi}_{1+k}(\mathbf{X})
  =
  \mathbf{Q}^{-1}
  \dfrac{\partial\boldsymbol{\psi}}{\partial X_k}(\mathbf{X}),
  \quad
  k=1,\ldots,N.
  \label{eq:phi-blocks}
\end{equation}
By construction,
$  \mathbb{E}\!\left[
    \mathbf{B}(\mathbf{X})
    \mathbf{B}^{\top}(\mathbf{X})
  \right]
  =
  \mathbf{Q}^{-1}
  \mathbf{G}_{\mathrm{GE}}
  \mathbf{Q}^{-\top}
  =
  \mathbf{Q}^{-1}(\mathbf{Q}\mathbf{Q}^\top)\mathbf{Q}^{-\top} 
  = 
  (\mathbf{Q}^{-1}\mathbf{Q})(\mathbf{Q}^\top\mathbf{Q}^{-\top}) 
  =
  \mathbf{I}_{n_\psi} \cdot \mathbf{I}_{n_\psi} 
  = 
  \mathbf{I}_{n_\psi}.
  \label{eq:ge-orthonormality}$
%
The first block
$\boldsymbol{\Phi}_1(\mathbf{X})$ is a
Sobolev-transformed representation of the original SDD
basis vector.
In general, $\boldsymbol{\Phi}_1(\mathbf{X})\ne\boldsymbol{\psi}(\mathbf{X})$,
because the derivative terms of
$\mathbf{G}_{\mathrm{GE}}$ transform the non-constant
basis directions.
We estimate the coefficient in this transformed
coordinate and then map them to the original $L^2$-orthonormal SDD coordinates for statistical moment evaluation.

\begin{remark}[Preservation of the constant component]
\label{rem:constant-retention}
Because $\psi_1\equiv1$, all partial derivatives of $\psi_1$ vanish.
The remaining SDD basis functions have zero mean and
are $L^2$-orthogonal to $\psi_1$.
Therefore, the first row and column of $\mathbf{G}_{\mathrm{GE}}$
are $(1,0,\ldots,0)^\top$ and
$(1,0,\ldots,0)$, respectively. The Cholesky
transformation thus preserves the constant basis direction and transforms only
the non-constant directions.
\end{remark}

Evaluating each transformed basis block at the $M$ training
samples introduced in~\eqref{eq:ge-training-data} gives
\begin{equation}
  \mathbf{A}_k
  :=
  \begin{bmatrix}
    \boldsymbol{\Phi}_k(\mathbf{x}^{(1)})^\top\\
    \boldsymbol{\Phi}_k(\mathbf{x}^{(2)})^\top\\
    \vdots\\
    \boldsymbol{\Phi}_k(\mathbf{x}^{(M)})^\top
  \end{bmatrix}
  \in
  \mathbb{R}^{M\times n_\psi},
  \qquad
  k=1,\ldots,1+N.
  \label{eq:block-design}
\end{equation}
Here, $\mathbf{A}_1$ contains the function-value evaluations, whereas
$\mathbf{A}_{k+1}$ contains the basis derivatives with respect to $X_k$ for
$k=1,\ldots,N$. Stacking these blocks gives 
\begin{equation}
  \Astack
  :=
  \begin{bmatrix}
    \mathbf{A}_1\\
    \mathbf{A}_2\\
    \vdots\\
    \mathbf{A}_{1+N}\\
  \end{bmatrix}
  \in
  \mathbb{R}^{M(1+N)\times n_\psi}.
  \label{eq:Astack}
\end{equation}
The corresponding target blocks are
\begin{align}
  \mathbf{b}_1
  &:=
  \bigl[
    y(\mathbf{x}^{(1)}),
    \ldots,
    y(\mathbf{x}^{(M)})
  \bigr]^\top
  \in
  \mathbb{R}^{M}.
  \\
  \mathbf{b}_{1+k}
  &:=
  \left[
    \frac{\partial y}{\partial x_k}
    \bigl(\mathbf{x}^{(1)}\bigr),
    \ldots,
    \frac{\partial y}{\partial x_k}
    \bigl(\mathbf{x}^{(M)}\bigr)
  \right]^\top
  \in
  \mathbb{R}^{M},
  \quad
  k=1,\ldots,N.
  \label{eq:b1k}
\end{align}
The stacked target vector is then
\begin{equation}
  \bstack
  :=
  \begin{bmatrix}
    \mathbf{b}_1\\
    \mathbf{b}_2\\
    \vdots\\
    \mathbf{b}_{1+N}
  \end{bmatrix}
  \in
  \mathbb{R}^{M(1+N)}.
  \label{eq:bstack}
\end{equation}

The resulting regression system is
\begin{equation}
  \Astack\mathbf{c}_{\mathrm{GE}}
  \approx
  \bstack,
\qquad
\mathbf{c}_{\mathrm{GE}}\in\mathbb{R}^{n_\psi},
\label{eq:stacked-system}
\end{equation}
where $\mathbf{c}_{\mathrm{GE}}$ denotes the coefficient vector
in the transformed Sobolev coordinates.
Each training sample contributes one
function-value equation and $N$ derivative
equations. Gradient augmentation can therefore improve
coefficient identifiability without additional training samples.

\subsection{Block-wise energy-balanced scaling}
\label{subsec:scaling}

The function-value and derivative design blocks can have substantially different Frobenius norms because of input units, physical scales, 
and basis-derivative magnitudes.
Without scaling, high-norm blocks can dominate the stacked regression
loss and reduce the accuracy of the function-value fit.

Define the block norms 
\begin{equation}
  F_f
  :=
  \|\mathbf{A}_1\|_F,
  \qquad
  F_{g_k}
  :=
  \|\mathbf{A}_{1+k}\|_F,
  \quad
  k=1,\ldots,N.
  \label{eq:block-norms}
\end{equation}
The corresponding block-wise scaling factors are
\begin{equation}
  s_f:=1,
  \qquad
  s_{g_k}
  :=
\begin{cases}
F_f/F_{g_k}, & F_{g_k}>\epsilon_s,\\
1, & F_{g_k}\le\epsilon_s,
\end{cases}
\quad
k=1,\ldots,N,
\label{eq:scale-factors}
\end{equation}
where $\epsilon_s>0$ prevents division by a near-zero derivative
norm.

We collect these factors in the diagonal row-weight matrix 
\begin{equation}
  \Wrow
  :=
  \operatorname{diag}\!\Bigl(
    \underbrace{
      s_f,\ldots,s_f,\;
      s_{g_1},\ldots,s_{g_1},\;
      \ldots,\;
      s_{g_N},\ldots,s_{g_N}
    }_{M(1+N)\text{ entries}}
  \Bigr)
  \in
  \mathbb{R}^{M(1+N) \times M(1+N)}.
  \label{eq:Wrow}
\end{equation}
Applying the same row weights to the design matrix and target
vector gives
\begin{equation}
  \Atil
  :=
  \Wrow\Astack,
  \qquad
  \btil
  :=
  \Wrow\bstack.
  \label{eq:balanced-system}
\end{equation}
For every derivative block satisfying $F_{g_k}>\epsilon_s$, the
scaling equalizes its Frobenius norm with that of the
function-value block:
\begin{equation}
  \left\|
    s_{g_k}\mathbf{A}_{1+k}
  \right\|_F
  =
  \|\mathbf{A}_1\|_F.
  \label{eq:balanced-block-norm}
\end{equation}

In the numerical examples, we compute the block norms and 
scaling factors once from the full training
design matrices and keep them fixed during grouped 
cross-validation and final model fitting. The scaling
uses only the sampled inputs and basis evaluations, not the
function-value or derivative targets.

\subsection{Ridge regression and grouped $K$-fold
cross-validation}
\label{subsec:ridge-cv}

In the limited-sample regime $M < n_\psi$, gradient
augmentation provides up to $M(N+1)$ equations, but the stacked
system may remain rank deficient or ill-conditioned. Ridge
regression stabilizes coefficient estimation by shrinking the
nonconstant coefficients. Unlike LASSO, ridge does not promote exact
sparsity and therefore allows gradient observations to inform all
retained basis directions.

We estimate the transformed-coordinate coefficients from the
balanced system by
\begin{equation}
  \hat{\mathbf{c}}_{\mathrm{GE}}(\alpha)
  :=
  \operatorname*{arg\,min}_{
    \mathbf{c}_{\mathrm{GE}}
    \in\mathbb{R}^{n_\psi}}
  \left\{
    \left\|
      \Atil\mathbf{c}_{\mathrm{GE}}-\btil
    \right\|_2^2
    +
    \alpha\,
    \mathbf{c}_{\mathrm{GE}}^\top
    \Dpen\,
    \mathbf{c}_{\mathrm{GE}}
  \right\},
  \qquad
  \alpha\ge0.
  \label{eq:ridge}
\end{equation}
The SDD basis already contains the constant function, so the
regression requires no separate intercept. We exclude the
constant coefficient from regularization using
$  \Dpen
  :=
  \operatorname{diag}(0,1,\ldots,1)
  \in\mathbb{R}^{n_\psi\times n_\psi}.
  \label{eq:penalty-matrix}$
The same constant-excluded penalty is used for the function-only
SDD baseline.

For $\alpha>0$, the solution is
\begin{equation}
  \hat{\mathbf{c}}_{\mathrm{GE}}(\alpha)
  =
  \left(
    \Atil^\top\Atil
    +
    \alpha\Dpen
  \right)^{-1}
  \Atil^\top\btil.
  \label{eq:ridge-solution}
\end{equation}
Although $\Dpen$ leaves the constant direction unpenalized, the
function-value block contains a nonzero constant column.
Therefore, the coefficient matrix in~\eqref{eq:ridge-solution}
is positive definite for $\alpha>0$.
We impose the ridge penalty in the transformed coordinates. The
transformed and original SDD coefficients satisfy $\mathbf{c}_{\mathrm{GE}}=\mathbf{Q}^{\top}\mathbf{c}_{\mathrm{SDD}}$. Therefore, 
\begin{equation}
  \mathbf{c}_{\mathrm{GE}}^\top
  \Dpen\,
  \mathbf{c}_{\mathrm{GE}}
  =
  \mathbf{c}_{\mathrm{SDD}}^\top
  \left(
    \mathbf{Q}\Dpen\mathbf{Q}^{\top}
  \right)
  \mathbf{c}_{\mathrm{SDD}}.
  \label{eq:sobolev-weighted-penalty}
\end{equation}
Because $\mathbf{Q}$ preserves the constant direction
(Remark~\ref{rem:constant-retention}), the transformed penalty
leaves the mean coefficient unregularized and applies a
Sobolev-weighted quadratic penalty only to the nonconstant SDD
coefficients.

We select the regularization parameter $\alpha$ using
grouped $K$-fold cross-validation.
All $N+1$ observations associated with the same training sample
$\mathbf{x}^{(i)}$ share group label
$i$.
Thus, the function value and its partial derivatives remain together in either the training or validation partition, preventing information leakage across folds. 

Let
$\mathcal{V}_r\subseteq\{1,\ldots,M\}$ denote the validation indices for fold $r$, and define
$\mathcal{T}_r := \{1,\ldots,M\}\setminus\mathcal{V}_r$. Under the block ordering in~\eqref{eq:Astack}, the rows of the
balanced stacked system associated with $\mathbf{x}^{(i)}$ occupy the rows 
 $ \mathcal{R}(i)
  :=
  \bigl\{\,
    i+(k-1)M
    \;:\;
    k=1,\ldots,1+N
  \,\bigr\}.
  \
  i=1,\ldots,M,
  \label{eq:row-index-map}$
Here, $k=1$ identifies the function-value row, and
$k=2,\ldots,1+N$ identify the derivative rows.

For fold $r$, denote the stacked training-row index set by
$\bigcup_{i\in\mathcal{T}_r}\mathcal{R}(i)$. The fold-specific training system is
\begin{equation}
  \widetilde{\mathbf{A}}_{\mathrm{tr}}^{(r)}
  :=
  \bigl[\Atil\bigr]_{
    \bigcup_{i\in\mathcal{T}_r}\mathcal{R}(i)}
  \in
  \mathbb{R}^{|\mathcal{T}_r|(1+N)\times n_\psi},
  \qquad
  \widetilde{\mathbf{b}}_{\mathrm{tr}}^{(r)}
  :=
  \bigl[\btil\bigr]_{
    \bigcup_{i\in\mathcal{T}_r}\mathcal{R}(i)}
  \in
  \mathbb{R}^{|\mathcal{T}_r|(1+N)},
  \label{eq:fold-train-blocks}
\end{equation}
where $[\,\cdot\,]_{\bigcup_{i\in\mathcal{T}_r}\mathcal{R}(i)}$
denotes the row-selection operator that extracts the rows indexed by
$\bigcup_{i\in\mathcal{T}_r}\mathcal{R}(i)$.
We then estimate
\begin{equation}
  \hat{\mathbf{c}}_{\mathrm{GE}}^{(r)}(\alpha)
  :=
  \operatorname*{arg\,min}_{\mathbf{c}\in\mathbb{R}^{n_\psi}}
  \left\{
    \left\|
      \widetilde{\mathbf{A}}_{\mathrm{tr}}^{(r)}
      \mathbf{c}
      -
      \widetilde{\mathbf{b}}_{\mathrm{tr}}^{(r)}
    \right\|_2^2
    +
    \alpha\,
    \mathbf{c}^\top\Dpen\,\mathbf{c}
  \right\}.
  \label{eq:ridge-fold}
\end{equation}

We evaluate validation error using only the held-out
function-value observations. Derivatives provide auxiliary
constraints during fitting but do not enter the validation loss,
which therefore measures function-prediction accuracy directly.
The validation blocks are
\begin{equation}
  \mathbf{A}_{1,\mathrm{val}}^{(r)}
  :=
  \bigl[\mathbf{A}_1\bigr]_{\mathcal{V}_r}
  \in
  \mathbb{R}^{|\mathcal{V}_r|\times n_\psi},
  \qquad
  \mathbf{b}_{1,\mathrm{val}}^{(r)}
  :=
  \bigl[\mathbf{b}_1\bigr]_{\mathcal{V}_r}
  \in
  \mathbb{R}^{|\mathcal{V}_r|},
  \label{eq:fold-val-blocks}
\end{equation}
where $\mathbf{A}_1$ and $\mathbf{b}_1$ are the
function-value design block in~\eqref{eq:block-design} and
target block in~\eqref{eq:b1k}, respectively.
The cross-validation loss is
\begin{equation}
  \mathcal{L}_{\mathrm{CV}}(\alpha)
  :=
  \frac{1}{K}
  \sum_{r=1}^{K}
  \frac{1}{|\mathcal{V}_r|}
  \left\|
    \mathbf{A}_{1,\mathrm{val}}^{(r)}
    \hat{\mathbf{c}}_{\mathrm{GE}}^{(r)}(\alpha)
    -
    \mathbf{b}_{1,\mathrm{val}}^{(r)}
  \right\|_2^2.
  \label{eq:cv-loss}
\end{equation}
We select $\alpha$ from a logarithmically spaced candidate grid
$\mathcal{G}$:
\begin{equation}
    \hat{\alpha}:=\operatorname*{arg\,min}_{\alpha\in\mathcal{G}}
    \mathcal{L}_{\mathrm{CV}}(\alpha).
\end{equation}
Finally, we refit~\eqref{eq:ridge} to the complete balanced
system using $\hat{\alpha}$.

\subsection{SDD coefficient and statistical moment
estimation}
\label{subsec:moments}

Ridge regression estimates the coefficients in the transformed
Sobolev coordinates. Moment estimation, however, requires
coefficients in the original $L^2$-orthonormal SDD coordinates.
From~\eqref{eq:phi-blocks}, 
$\boldsymbol{\Phi}_1(\mathbf{X})
  =
  \mathbf{Q}^{-1}
  \boldsymbol{\psi}(\mathbf{X}).
  \label{eq:phi1-relation}$
The fitted GE-SDD surrogate is therefore $
  \hat{y}(\mathbf{X})=
  \boldsymbol{\Phi}_1(\mathbf{X})^\top
  \hat{\mathbf{c}}_{\mathrm{GE}}
  =
  \boldsymbol{\psi}(\mathbf{X})^\top
  \mathbf{Q}^{-\top}
  \hat{\mathbf{c}}_{\mathrm{GE}}.
$ Accordingly, we estimate the coefficients in the original SDD
coordinates as
\begin{equation}
    \hat{\mathbf{c}}_{\mathrm{SDD}}
    :=
    \mathbf{Q}^{-\top}
    \hat{\mathbf{c}}_{\mathrm{GE}}.
  \label{eq:coef-estimation}
\end{equation}
The SDD surrogate then becomes
\begin{equation}
  \hat{y}(\mathbf{X})
  =
  \boldsymbol{\psi}(\mathbf{X})^\top
  \hat{\mathbf{c}}_{\mathrm{SDD}}.
  \label{eq:sdd-estimated-surrogate}
\end{equation}
Because $\boldsymbol{\psi}$ is $L^2$-orthonormal with respect to
the input probability measure, the mean and variance follow
directly from the estimated coefficients:
\begin{equation}
  \hat{\mu}_Y
  =
  \hat{c}_{\mathrm{SDD},1},
  \qquad
  \hat{\sigma}_Y^2
  =
  \sum_{j=2}^{n_\psi}
  \hat{c}_{\mathrm{SDD},j}^{\,2}.
  \label{eq:estimated-moments}
\end{equation}
Since $\psi_1\equiv1$, the coefficient
$\hat{c}_{\mathrm{SDD},1}=\hat{c}_\emptyset$ estimates the ANOVA
constant $y_\emptyset=\mathbb{E}[Y]$. Thus, GE-SDD provides
moment estimates without additional Monte Carlo sampling.

The same SDD coefficients also represent the fitted
gradient:
\begin{align}
  \dfrac{\partial}{\partial{X_k}}\hat{y}(\mathbf{X})
  &=
  \boldsymbol{\Phi}_{1+k}(\mathbf{X})^\top
  \hat{\mathbf{c}}_{\mathrm{GE}}
  =
  \left(
    \dfrac{\partial}{\partial{X_k}}\boldsymbol{\psi}(\mathbf{X})
  \right)^\top
  \hat{\mathbf{c}}_{\mathrm{SDD}},
  \qquad
  k=1,\ldots,N.
  \label{eq:gradient-consistency}
\end{align}
The reconstructed SDD coefficient vector 
$\hat{\mathbf{c}}_{\mathrm{SDD}}$ therefore simultaneously and 
consistently represents both the fitted response and its 
 first-order gradient derivatives.

\subsection{Computational workflow and cost}
\label{subsec:computational-procedure}

This subsection outlines the computational procedure and the
associated cost of constructing the $S$-variate GE-SDD surrogate.
Figure~\ref{fig:flowchart} presents the overall workflow, and
Algorithm~\ref{alg:gesdd} states the complete procedure at the end
of this subsection.

The procedure begins with the input-output-gradient dataset
$\{\mathbf{x}^{(i)},y(\mathbf{x}^{(i)}),
\nabla y(\mathbf{x}^{(i)})\}_{i=1}^{M}$. We first construct the
$L^2$-orthonormal SDD basis $\boldsymbol{\psi}(\mathbf{X})$ from the
prescribed spline parameters through the whitening transformation
described in Section~\ref{subsec:sdd} and
Appendix~\ref{app:basis-construction}. We then assemble the extended
basis matrix $\mathbf{A}(\mathbf{X})$ from the basis and its partial
derivatives and factorize the derivative-augmented spline moment
matrix $\mathbf{G}_{\mathrm{GE}}=\mathbf{Q}\mathbf{Q}^\top$, which
maps the function value and its gradients into the common Sobolev
coordinates $\boldsymbol{\Phi}_1,\ldots,\boldsymbol{\Phi}_{1+N}$
(Section~\ref{subsec:stacked}).

Because the function-value and derivative blocks generally differ in
Frobenius norm, we introduce the block-wise scaling factors
$s_{g_k}=\|\mathbf{A}_1\|_F/\|\mathbf{A}_{1+k}\|_F$ to form the
energy-balanced system $(\Atil,\btil)$, which prevents a high-norm
block from dominating the stacked regression loss
(Section~\ref{subsec:scaling}). We then estimate the coefficients by
constant-excluded ridge regression, selecting $\hat{\alpha}$ by
grouped $K$-fold cross-validation that keeps all $1+N$ observations of
each training point in the same fold
(Section~\ref{subsec:ridge-cv}). Finally, we map the solution back to
the $L^2$-orthonormal SDD coordinates
($\hat{\mathbf{c}}_{\mathrm{SDD}}=\mathbf{Q}^{-\top}
\hat{\mathbf{c}}_{\mathrm{GE}}$), preserving closed-form mean and
variance estimates (Section~\ref{subsec:moments}). The cost is
governed by $n_\psi=L_{S,\mathbf{p},\boldsymbol{\Xi}}$
in~\eqref{eq:basis-count} and the $M(1+N)$ stacked equations, with
detailed analysis deferred to Section~\ref{subsec:discussion}.
\begin{algorithm}[H]
\caption{Gradient-enhanced SDD (GE-SDD)}
\label{alg:gesdd}
\begin{algorithmic}[1]

\Require
Training data
$\bigl\{\bigl(\mathbf{x}^{(i)},
y(\mathbf{x}^{(i)}),
\nabla y(\mathbf{x}^{(i)})\bigr)\bigr\}_{i=1}^{M}$;
spline parameters
$\{p_k,n_k,\boldsymbol{\xi}_k\}_{k=1}^{N}$;
interaction order $S$;
number of CV folds $K$;
regularization grid $\mathcal{G}$

\Ensure
$\hat{\mathbf{c}}_{\mathrm{GE}}$,
$\hat{\mathbf{c}}_{\mathrm{SDD}}$,
$\hat{\mu}_Y$,
$\hat{\sigma}_Y^2$

\Statex \textit{Step 1: Orthonormal basis construction}
\For{$k=1,\ldots,N$}
  \State Assemble the raw basis vector
  $\mathbf{P}_k(X_k)$ (Appendix~\ref{app:basis-construction}).
  \State Compute
  $\mathbf{G}_k
  =\mathbb{E}\bigl[\mathbf{P}_k\mathbf{P}_k^\top\bigr]$
  and factorize
  $\mathbf{G}_k=\mathbf{Q}_k\mathbf{Q}_k^\top$.
  \State Form
  $\boldsymbol{\psi}_k(X_k)
  =\mathbf{Q}_k^{-1}\mathbf{P}_k(X_k)$.
\EndFor
\State Construct the retained tensor-product basis
$\Psi^u_{\mathbf{i}_u}
=\prod_{k\in u}\psi^k_{i_k}$,
$\mathbf{i}_u\in\bar{\mathcal{I}}_{u,\mathbf{n}_u}$,
$1\le|u|\le S$,
\Statex \hspace{\algorithmicindent}
and collect the $n_\psi$ functions into
$\boldsymbol{\psi}(\mathbf{X})$
(Section~\ref{subsec:sdd}).

\Statex \textit{Step 2: Sobolev transformation}
\State Assemble
$\mathbf{A}(\mathbf{X})
=\bigl[\,\boldsymbol{\psi}~~
\partial\boldsymbol{\psi}/\partial X_1~~
\cdots~~
\partial\boldsymbol{\psi}/\partial X_N\,\bigr]$.
\State Compute
$\mathbf{G}_{\mathrm{GE}}
=\mathbb{E}\bigl[\mathbf{A}\mathbf{A}^\top\bigr]$
and factorize
$\mathbf{G}_{\mathrm{GE}}=\mathbf{Q}\mathbf{Q}^\top$
(Section~\ref{subsec:stacked}).
\State Form
$\boldsymbol{\Phi}_1=\mathbf{Q}^{-1}\boldsymbol{\psi}$
and
$\boldsymbol{\Phi}_{1+k}
=\mathbf{Q}^{-1}\,
\partial\boldsymbol{\psi}/\partial X_k$,
$k=1,\ldots,N$.

\Statex \textit{Step 3: Stacking and energy-balanced scaling}
\State Evaluate
$\boldsymbol{\Phi}_1,\ldots,\boldsymbol{\Phi}_{1+N}$
at the $M$ training points; stack into
$\Astack$ and $\bstack$.
\State Compute
$s_f=1$ and
$s_{g_k}=\|\mathbf{A}_1\|_F/\|\mathbf{A}_{1+k}\|_F$;
collect them into $\Wrow$.
\State Form the balanced system
$\Atil=\Wrow\Astack$ and
$\btil=\Wrow\bstack$
(Section~\ref{subsec:scaling}).

\Statex \textit{Step 4: Ridge regression with grouped cross-validation}
\State Assign all $1+N$ rows of each training point
$\mathbf{x}^{(i)}$ to group $i$
(Section~\ref{subsec:ridge-cv}).
\For{each $\alpha\in\mathcal{G}$}
  \For{$r=1,\ldots,K$}
    \State Fit the constant-excluded ridge estimator
    $\hat{\mathbf{c}}_{\mathrm{GE}}^{(r)}(\alpha)$
    on the training groups of fold $r$.
  \EndFor
  \State Evaluate the function-value validation loss
  $\mathcal{L}_{\mathrm{CV}}(\alpha)$
  over the held-out function observations.
\EndFor
\State Select
$\hat{\alpha}
=\arg\min_{\alpha\in\mathcal{G}}
\mathcal{L}_{\mathrm{CV}}(\alpha)$.
\State Refit on the complete balanced system
$(\Atil,\btil)$ to obtain
$\hat{\mathbf{c}}_{\mathrm{GE}}$.

\Statex \textit{Step 5: Coefficient and moment estimation}
\State Transform back:
$\hat{\mathbf{c}}_{\mathrm{SDD}}
=\mathbf{Q}^{-\top}\hat{\mathbf{c}}_{\mathrm{GE}}$
(Section~\ref{subsec:moments}).
\State Compute
$\hat{\mu}_Y=\hat{c}_{\mathrm{SDD},1}$
and
$\hat{\sigma}_Y^2
=\sum_{j=2}^{n_\psi}\hat{c}_{\mathrm{SDD},j}^{\,2}$.

\end{algorithmic}
\end{algorithm} 

\section{Numerical examples}
\label{sec:numerical}

We evaluate the proposed GE-SDD on three
benchmarks of increasing input dimension and
complexity: a two-dimensional continuous exponential
function in Section~\ref{subsec:ex1}, a two-degree-of-freedom linear dynamical
system with three uncertain parameters in Section~\ref{subsec:ex2}, and a
30-dimensional compliance problem for a 25-bar transmission-tower-type space truss in Section~\ref{subsec:ex3}.
Unless stated otherwise, all experiments follow the common setup below.
Tables report
arithmetic means over 20 independent replications,
while box plots summarize the distributions and
medians. We compare GE-SDD with:
\begin{itemize}
  \item \textbf{SDD}: conventional SDD trained on function
        values, using the same orthonormal basis,
        constant-excluded ridge estimator, regularization grid,
        and $K$-fold cross-validation procedure as GE-SDD.
  \item \textbf{Kriging}: Gaussian process regression trained
        on function values~\cite{Bouhlel2019}. The model uses a constant trend, affinely scaled inputs, standardized outputs, and anisotropic length scales
        estimated by maximum likelihood.
        Example~1 uses a Mat\'ern~$3/2$ kernel with a nugget of
        $10^{-10}$ and five random restarts. Examples~2 and~3
        use an anisotropic squared-exponential kernel with a
        nugget of $10^{-8}$. 
  \item \textbf{GE-Kriging}: Gradient-enhanced Kriging based on
        the exact joint covariance of function values and gradients for the
        corresponding kernel~\cite{Morris1993,Bouhlel2019}. GE-Kriging uses the
        same $M$ training samples and $M\times N$ derivative
        observations as GE-SDD. Each example specifies the hyperparameter-selection
        strategy.
  \item \textbf{RBF}: A Gaussian radial basis function surrogate,
        $\phi(\mathbf{x},\mathbf{x}')
        =\exp(-\|\mathbf{x}-\mathbf{x}'\|_2^2/d_0^2)$,
        trained on function values and included only in Example~1.
        We select the shape parameter $d_0$ by $K$-fold
        cross-validation over a logarithmically spaced grid. The model uses the same input and output scaling as Kriging, no
        polynomial tail, and a regularization coefficient of
        $10^{-10}$.
\end{itemize}

We measure surrogate accuracy using the normalized
root-mean-square error
$  \mathrm{NRMSE}
  :=
  \frac{\mathrm{RMSE}}
       {\sigma_{Y,\mathrm{MCS}}}
  \times 100\%,
  \label{eq:nrmse}$
where the RMSE is evaluated on a reference MCS sample
of size $M_{\mathrm{MCS}}=10^6$, or the fixed subsample
specified for an individual example. The quantity  $\sigma_{Y,\mathrm{MCS}}$ is
the corresponding MCS estimate of the output standard deviation. We also report the coefficient of determination $R^2$ and the mean
absolute error (MAE) where
relevant.

\subsection{Example 1: Two-dimensional continuous
exponential function}
\label{subsec:ex1}

The first benchmark is a two-dimensional continuous
exponential function with a localized peak at the origin: 
\begin{equation}
  y(x_1,x_2)
  =
  \exp\!\bigl(-2|x_1|-2|x_2|\bigr),
  \qquad
  x_1,x_2\in[-1,1].
  \label{eq:cont-exp}
\end{equation}
The inputs $X_1$ and $X_2$ are independent and  uniformly distributed on $[-1,1]$.
The response decays rapidly from the origin and is continuous but nondifferentiable along the coordinate axes. These localized, nonsmooth features challenge surrogates trained only on function values.
Away from the coordinate axes, the partial derivatives are
\begin{equation}
  \frac{\partial y}{\partial x_k}
  =
  -2\,\operatorname{sign}(x_k)
  \exp\!\bigl(-2|x_1|-2|x_2|\bigr),
  \qquad
  k=1,2,\quad x_k\ne0.
  \label{eq:cont-exp-grad}
\end{equation}
The classical derivative is undefined at $x_k=0$. 
We set
$\operatorname{sign}(0)=0$, corresponding to the
symmetric subgradient convention at the cusp.
Because the training samples follow continuous uniform distributions, the probability
of sampling $x_k=0$ exactly is zero;
the adopted convention therefore does not affect the
reported random-sampling results.

\paragraph{Spline discretization.}
SDD and GE-SDD use the same
cubic B-spline basis with $p_k=3$ for $k=1,2$.
Each input coordinate uses the open knot vector
$  \boldsymbol{\xi}_k
  =
  \bigl(
    \underbrace{-1,\ldots,-1}_{4},\,
    -\tfrac{1}{2},\,
    \underbrace{0,\,0,\,0}_{3},\,
    \tfrac{1}{2},\,
    \underbrace{1,\ldots,1}_{4}
  \bigr),
  \
  k=1,2.
  \label{eq:example1-knots}$
The interior knots are 
$\{-\tfrac{1}{2},0,\tfrac{1}{2}\}$ with multiplicity three at the
origin. 
This repeated knot reduces the interelement
continuity of the cubic spline space from $C^{2}$ to
$C^{3-3}=C^{0}$ at $x_k=0$, so the univariate spline
space can represent a kink exactly at the known location
of the nonsmooth ridge of~\eqref{eq:cont-exp}. Both SDD and GE-SDD methods use the same knot vector, so their comparison isolates the effect of gradient observations from that of the
approximation space.
The knot vector in~\eqref{eq:example1-knots} generates
$n_k=9$ univariate basis functions per coordinate.

With interaction order $S=2$, the number of retained
basis functions is
$  n_\psi
  =
  1+2(9-1)+(9-1)^2
  =
  81.
  \label{eq:example1-basis-count}$
Both surrogates use the same $M=81$ independently
drawn training samples. The function-only SDD design matrix
is therefore square, with 81 equations and $81$ coefficients. GE-SDD adds two derivative equations per sample, producing $3M=243$ stacked equations for the same
$81$ coefficients.

\begin{figure}
  \centering
  \includegraphics[width=.9\textwidth]
  {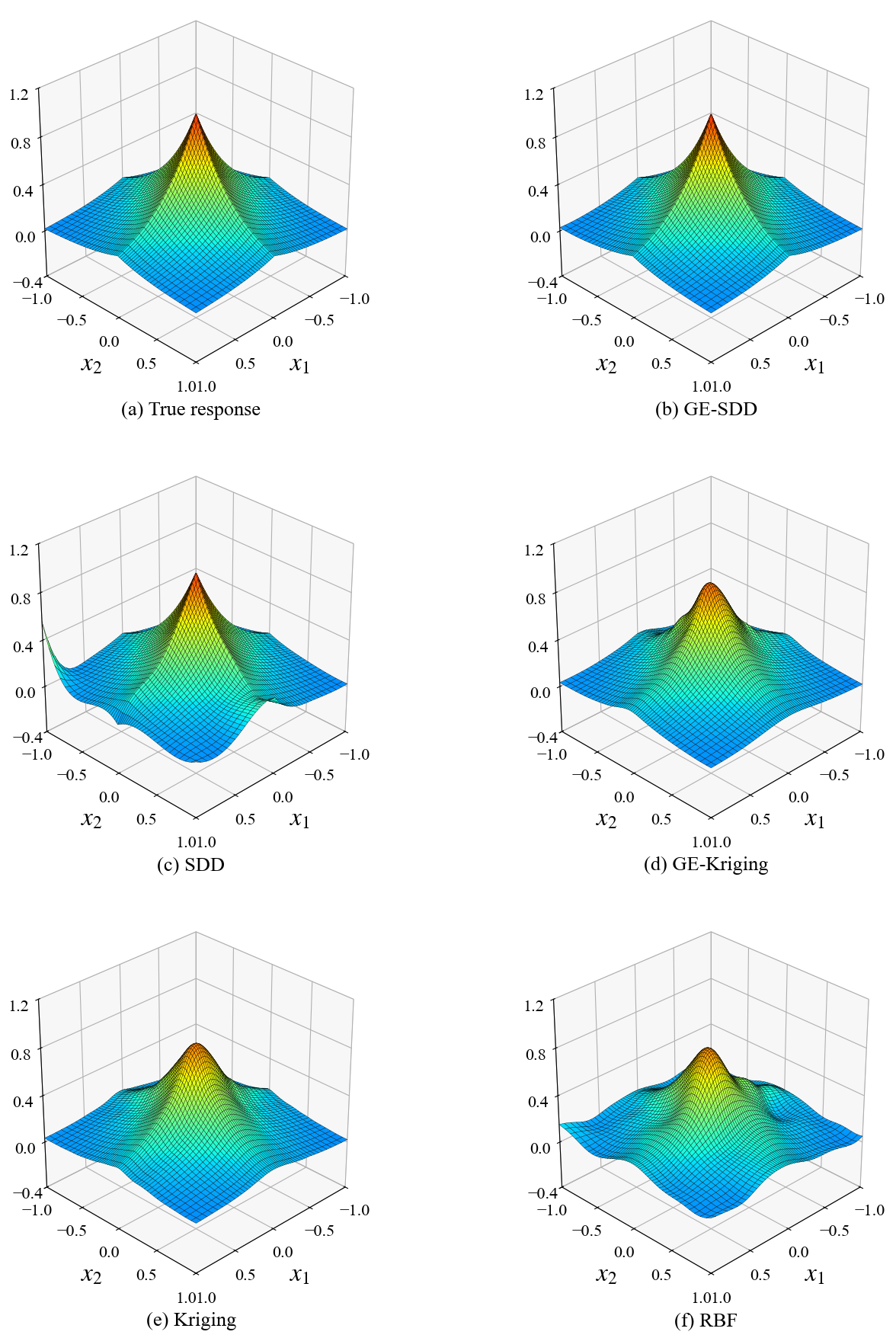}
  \caption{True and surrogate-predicted response surfaces for
  the two-dimensional continuous exponential function using the
  representative training-data realization generated with seed 
  42:
          \textup{(a)} true response;
          \textup{(b)} GE-SDD (NRMSE $0.320\%$);
          \textup{(c)} standard SDD ($15.468\%$);
          \textup{(d)} GE-Kriging ($7.300\%$);
          \textup{(e)} Kriging ($7.735\%$); and
          \textup{(f)} RBF ($18.710\%$).}
  \label{fig:exp-surfaces}
\end{figure}

\clearpage

\begin{figure}
  \centering
  \includegraphics[width=.9\textwidth]
  {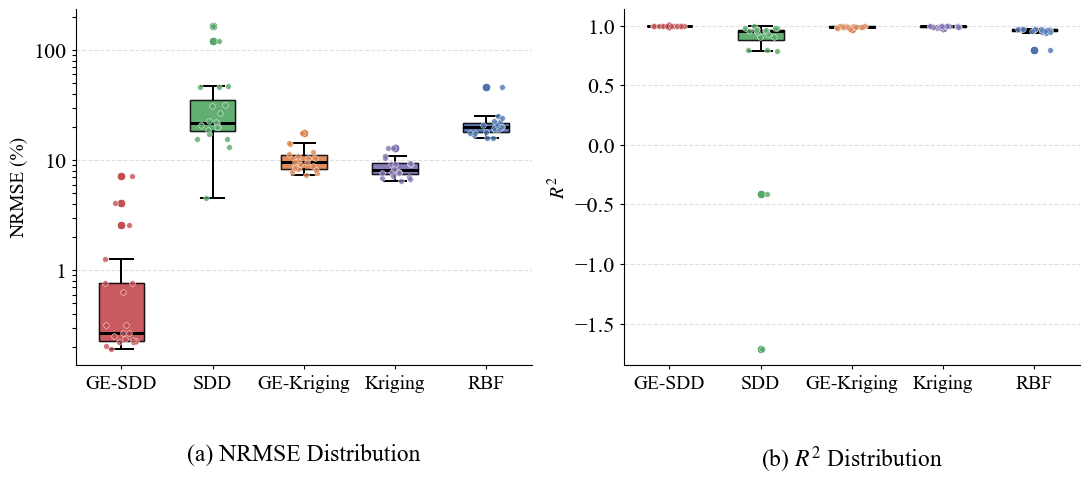}
    \caption{Distributions of pointwise prediction-accuracy metrics
    over 20 independent training-data realizations for the
    two-dimensional continuous exponential function
    ($M=81$, $p_k=3$, $n_k=9$, and $S=2$):
    \textup{(a)} NRMSE on a logarithmic scale and
    \textup{(b)} $R^2$.
    The overlaid points represent the results from individual
    realizations, and the horizontal line within each box denotes
    the median.}
  \label{fig:exp-boxplot}
\end{figure}

\begin{table}[ht]
\centering
    \caption{Mean pointwise prediction-accuracy metrics over
    20 independently generated training sets for the
    two-dimensional continuous exponential function
    ($M=81$, $p_k=3$, $n_k=9$, and $S=2$).
    NRMSE is reported as a percentage.}
\label{tab:example1-results}
\begin{tabular*}{\textwidth}
{@{\extracolsep{\fill}}lcc@{}}
\toprule
Method & NRMSE (\%) & $R^2$ \\
\midrule
GE-SDD     &  1.022 & 0.9996 \\
SDD        & 36.326 & 0.7298 \\
GE-Kriging & 10.238 & 0.9891 \\
Kriging    &  8.731 & 0.9922 \\
RBF        & 21.138 & 0.9526 \\
\bottomrule
\end{tabular*}
\end{table}

Table~\ref{tab:example1-results} reports the mean
prediction accuracy over the 20 replications.
GE-SDD achieves the lowest mean NRMSE of $1.022\%$ and the
highest mean $R^2$ of $0.9996$. Its mean NRMSE is more than eight times lower than that of Kriging ($8.731\%$),
the most accurate baseline, and more than
$35$ times lower than that of  standard SDD ($36.326\%$).
Function-only Kriging achieves slightly higher accuracy
GE-Kriging ($8.731\%$ versus $10.238\%$).
One
plausible explanation is the mismatch between the Mat\'ern~$3/2$
joint covariance model and the nonsmooth response. The true
gradient in~\eqref{eq:cont-exp-grad} changes sign discontinuously
across the coordinate axes, whereas the covariance model imposes
a smooth joint function--gradient structure. Gradient
observations on opposite sides of a cusp may therefore conflict
with the assumed covariance. In contrast, the GE-SDD basis
admits $C^0$ continuity at the known cusp locations, allowing
the model to use the gradient sign changes directly.

Figure~\ref{fig:exp-boxplot} summarizes the variability of the metrics across replications.
The median NRMSEs are $0.269\%$ for GE-SDD,
$21.772\%$ for SDD, $9.683\%$ for GE-Kriging, $8.086\%$ for
Kriging, and $19.945\%$ for RBF. GE-SDD achieves an NRMSE below 1\% in 15 of the 20 replications. The largest GE-SDD NRMSE of $7.208\%$ remains below the median Kriging NRMSE. Standard SDD produces NRMSEs ranging from $4.5\%$ to $166.6\%$ and negative $R^2$ values in 2 replications. These results show that the square function-only system is sensitive to the placement of training samples near the localized peak. The derivative constraints substantially reduce this sensitivity.

The response surfaces in Figure~\ref{fig:exp-surfaces}
support the quantitative results.
GE-SDD captures the peak amplitude and rapid decay in Figure~\ref{fig:exp-surfaces}(b), whereas standard SDD produces local artifacts near the peak and in sparsely sampled regions in Figure~\ref{fig:exp-surfaces}(c).
GE-Kriging and Kriging reproduce the reference response away
from the origin but smooth the cusps. RBF shows the
largest deviations near the peak and the domain boundary.

\subsection{Example 2: Two-degree-of-freedom linear
dynamical system}
\label{subsec:ex2}

\begin{figure}
  \centering
  \includegraphics[width=.7\columnwidth]
  {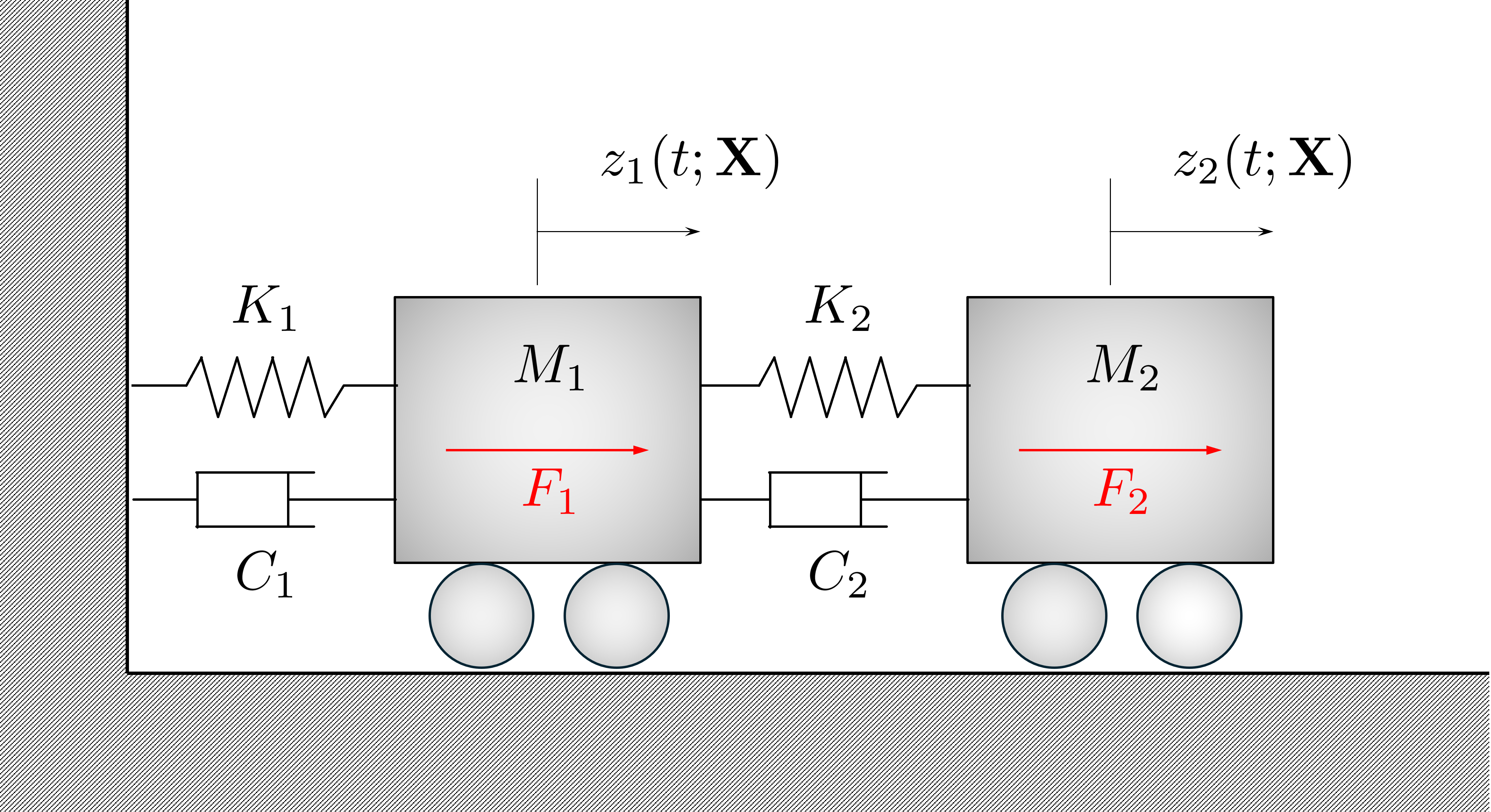}
  \caption{Two-degree-of-freedom linear dynamical system with
masses $M_1$ and $M_2$, stiffnesses $K_1$ and $K_2$, damping
coefficients $C_1$ and $C_2$, forces $F_1$ and $F_2$, and
displacements $z_1(t;\mathbf{X})$ and $z_2(t;\mathbf{X})$.}
  \label{fig:2dof-schematic}
\end{figure}

The second benchmark is a two-degree-of-freedom (2-DOF)
linear dynamical system with three uncertain structural
parameters representing mass, damping, and stiffness. Figure~\ref{fig:2dof-schematic} illustrates the system configuration.
The spring--damper
pair $(K_1,C_1)$ connects mass $M_1$ to the fixed support, while
$(K_2,C_2)$ connects mass $M_2$ to $M_1$. The forces $F_1(t)$ and
$F_2(t)$ act in the positive displacement direction, and
$Z_1(t;\mathbf{X})$ and $Z_2(t;\mathbf{X})$ denote the
displacements.

For excitation frequency $f$ in Hz, define the angular frequency
$  \omega:=2\pi f.
  \label{eq:angular-frequency}$
The frequency-domain equilibrium equation is
\begin{equation}
  \mathbf{D}(f;\mathbf{X})\,
  \mathbf{Z}(f;\mathbf{X})
  =
  \widehat{\mathbf{F}},
  \label{eq:2dof-eom}
\end{equation}
where
\begin{equation}
  \mathbf{D}(f;\mathbf{X})
  :=
  -\omega^2
  \begin{bmatrix}
    M_1 & 0\\
    0 & M_2
  \end{bmatrix}
  +
  \mathrm{i}\omega
  \begin{bmatrix}
    C_1+C_2 & -C_2\\
    -C_2 & C_2
  \end{bmatrix}
  +
  \begin{bmatrix}
    K_1+K_2 & -K_2\\
    -K_2 & K_2
  \end{bmatrix},
  \label{eq:2dof-dynamic-stiffness}
\end{equation}
and
\begin{equation}
  \mathbf{Z}(f;\mathbf{X})
  :=
  \begin{bmatrix}
    Z_1(f;\mathbf{X})\\
    Z_2(f;\mathbf{X})
  \end{bmatrix},
  \qquad
  \widehat{\mathbf{F}}
  :=
  \begin{bmatrix}
    \widehat F_1\\
    \widehat F_2
  \end{bmatrix}.
  \label{eq:2dof-response-force}
\end{equation}

The uncertain input vector is
\begin{equation}
  \mathbf{X}
  :=
  (X_M,X_C,X_K)^\top
  \in\mathbb{R}^3,
  \label{eq:2dof-input}
\end{equation}
where $X_M$, $X_C$, and $X_K$ are independent
$U(-3,3)$ random variables. The physical parameters are
\begin{equation}
  M_1=M_2
  =
  (1+0.05X_M)\;\mathrm{kg},
  \qquad
  C_1=C_2
  =
  (1+0.05X_C)\;\mathrm{N{\cdot}s/m},
  \qquad
  K_1=K_2
  =
  15000(1+0.05X_K)\;\mathrm{N/m}.
  \label{eq:2dof-params}
\end{equation}

We apply a unit harmonic force to $M_1$, such that
 $ \widehat{\mathbf F}
  =
  (1,0)^\top.
  \label{eq:2dof-numerical-force}$
%
At each fixed frequency, the surrogate approximates the response
magnitude
\begin{equation}
  y_f(\mathbf{X})
  :=
  \left|Z_1(f;\mathbf{X})\right|.
  \label{eq:2dof-qoi}
\end{equation}
We construct a separate surrogate at each of 100 equally spaced
frequencies
$  f_j\in[10,35]\;\mathrm{Hz},
  \ j=1,\ldots,100.$
The estimated coefficients provide the frequency-dependent
standard deviation
\begin{equation}
  \sigma_{|Z_1|}(f_j)
  :=
  \sqrt{
    \mathbb{V}\mathrm{ar}
    \left[
      |Z_1(f_j;\mathbf X)|
    \right]
  }.
  \label{eq:2dof-std-curve}
\end{equation}

\paragraph{Analytic response sensitivities.}
We compute the GE-SDD gradients by differentiating the complex
equilibrium equation.
For
$q\in\{X_M,X_C,X_K\}$, differentiation of
\eqref{eq:2dof-eom} gives
\begin{equation}
  \frac{\partial\mathbf Z}{\partial q}
  =
  -\mathbf D^{-1}
  \frac{\partial\mathbf D}{\partial q}
  \mathbf Z,
  \label{eq:2dof-response-sensitivity}
\end{equation}
because $\widehat{\mathbf{F}}$ does not depend on the uncertain
parameters. 
The required derivatives are
\begin{equation}
  \frac{\partial\mathbf D}{\partial X_M}
  =
  -\omega^2
  \left(
    0.05\,\mathbf I_2
  \right),
  \qquad
  \frac{\partial\mathbf D}{\partial X_C}
  =
  \mathrm{i}\omega
  \begin{bmatrix}
    0.10 & -0.05\\
    -0.05 & 0.05
  \end{bmatrix},
  \qquad
  \frac{\partial\mathbf D}{\partial X_K}
  =
  \begin{bmatrix}
    1500 & -750\\
    -750 & 750
  \end{bmatrix}.
  \label{eq:2dof-dD}
\end{equation}
For $|Z_1|>0$, the response-magnitude derivative is
\begin{equation}
  \frac{\partial y_f}{\partial q}
  =
  \frac{
    \operatorname{Re}\!\left[
      \overline{Z_1}
      \dfrac{\partial Z_1}{\partial q}
    \right]
  }{
    |Z_1|
  },
  \qquad
  q\in\{X_M,X_C,X_K\}.
  \label{eq:2dof-magnitude-sensitivity}
\end{equation}
We set the derivative to zero when $|Z_1|\le10^{-15}$ to avoid
division by a near-zero magnitude. Equations
\eqref{eq:2dof-response-sensitivity}--%
\eqref{eq:2dof-magnitude-sensitivity} define a direct analytical
sensitivity calculation rather than an adjoint method.

\paragraph{Spline discretization.}
SDD and GE-SDD  use $p_k=2$, $n_k=18$, and maximum interaction
order $S=2$ in all 3 input dimensions. The retained basis
contains
\begin{equation}
  n_\psi
  =
  1+3(18-1)+{3\choose2}(18-1)^2
  =
  919
  \label{eq:2dof-basis-count}
\end{equation}
functions. Each frequency-wise surrogate uses $M=919$
training samples.

\paragraph{Baseline configuration and evaluation.}
We train an independent Kriging model at each frequency using all
$M=919$ training samples and estimate its amplitude and
anisotropic length scales by maximum likelihood. GE-Kriging uses
the same $MN$ analytical derivative observations through the
joint function--gradient covariance. Repeated maximum-likelihood
optimization of the $(N+1)M\times(N+1)M$ covariance matrix is
impractical across 100 frequencies. We therefore select the
GE-Kriging hyperparameters by training-sample-grouped $K$-fold
cross-validation over a logarithmically spaced grid.

We estimate the Kriging and GE-Kriging standard-deviation curves
from predictions on a fixed $5\times10^4$-sample subset of the
Monte Carlo reference set. SDD and GE-SDD provide the same
statistics directly from their coefficients using
\eqref{eq:estimated-moments}.

\begin{figure}
\centering
\includegraphics[width=.75\columnwidth]
{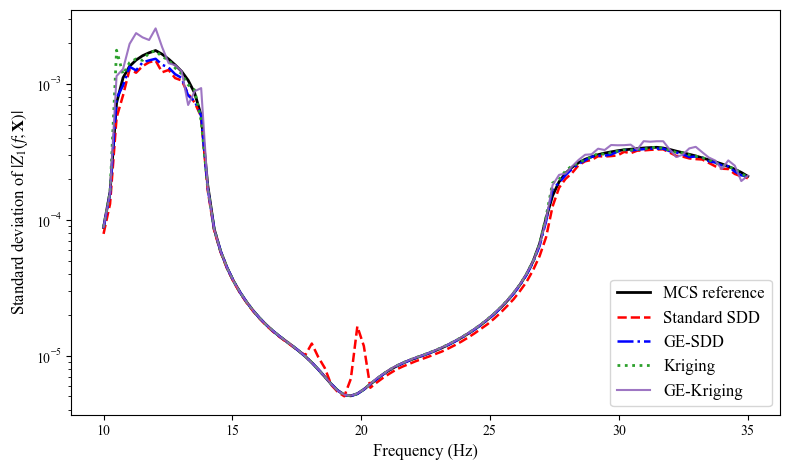}
\caption{Frequency-dependent standard deviation of
$|Z_1(f;\mathbf{X})|$ over $10$--$35\,\mathrm{Hz}$ for the
training set generated with seed 42. The Monte Carlo reference
is compared with SDD, GE-SDD, Kriging, and GE-Kriging
($M=919$, $n_\psi=919$, $S=2$, and
$M_{\mathrm{MCS}}=10^6$).}
\label{fig:2dof-freq}
\end{figure}

\begin{table}
\caption{Accuracy of the frequency-dependent standard deviation
of $|Z_1(f;\mathbf{X})|$ over 20 independently generated
training sets ($M=919$, $n_\psi=919$, and $S=2$). The metrics
are the mean absolute error, mean relative error across
frequencies, and coefficient of determination. The Monte Carlo
reference uses $M_{\mathrm{MCS}}=10^6$ samples.}
\label{tab:2dof}
\begin{tabular*}{\textwidth}{@{\extracolsep{\fill}}lrrr@{}}
\toprule
Method & MAE ($\times10^{-5}$) & Rel.\ err.\ (\%) & $R^2$ \\
\midrule
GE-SDD     & 2.435 &  2.45 & 0.977 \\
SDD        & 4.978 & 15.95 & 0.387 \\
Kriging    & 1.151 &  1.49 & 0.994 \\
GE-Kriging & 4.651 &  5.97 & 0.870 \\
\bottomrule
\end{tabular*}
\end{table}

This benchmark tests a frequency-dependent response with strong
resonance behavior. Near resonance, small parameter changes shift
both the magnitude and location of the response peaks, producing
a strongly nonlinear input--output map despite the linear
equations of motion. Figure~\ref{fig:2dof-freq} shows that all
4 surrogates agree closely with the reference away from
resonance. The main differences occur near the first resonance
at approximately $12\,\mathrm{Hz}$. Standard SDD underestimates
the peak, consistent with shrinkage of the nonconstant
coefficients under the regularization selected for the square
function-only system. GE-Kriging shows the largest local
deviations, while GE-SDD and Kriging track the reference closely.

Table~\ref{tab:2dof} reports the mean accuracy over 20
replications. Kriging achieves the lowest mean MAE of
$1.151\times10^{-5}$. GE-SDD follows with
$2.435\times10^{-5}$, GE-Kriging with
$4.651\times10^{-5}$, and standard SDD with
$4.978\times10^{-5}$.

GE-SDD is more accurate and less variable than standard SDD and
GE-Kriging in this benchmark. GE-SDD produces a lower MAE than
both methods in all 20 replications and reduces the mean relative
error of standard SDD from $15.95\%$ to $2.45\%$. The coefficient
of variation of the GE-SDD MAE is $9.4\%$, compared with
$20.1\%$ for Kriging, $34.2\%$ for GE-Kriging, and $68.3\%$ for
standard SDD.

The same gradient observations provide less benefit to
GE-Kriging. GE-Kriging is less accurate than function-only
Kriging on average, and its least accurate replication produces
a curve $R^2$ of $0.118$. This result indicates the difficulty of
selecting joint-covariance hyperparameters near resonance.

Function-only Kriging remains the most accurate method for this
low-dimensional benchmark and achieves approximately half the
GE-SDD MAE in every replication. This advantage is consistent
with frequency-specific maximum-likelihood adaptation of the
anisotropic length scales. Near resonance, the estimated length
scales decrease by more than an order of magnitude in the mass
and stiffness directions while remaining large in the damping
direction. This behavior follows the resonance relation
$\omega_{\mathrm{res}}^2\propto K/M$. The fixed spline basis used
by SDD and GE-SDD does not adapt with frequency, so the remaining
GE-SDD error concentrates near the first resonance. This result
motivates the gradient-guided adaptive knot placement discussed
in Section~\ref{sec:conclusions}.

\subsection{Example 3: 30-dimensional 25-bar space truss
compliance problem}
\label{subsec:ex3}

\begin{figure}
  \centering
  \includegraphics[width=.5\columnwidth]
  {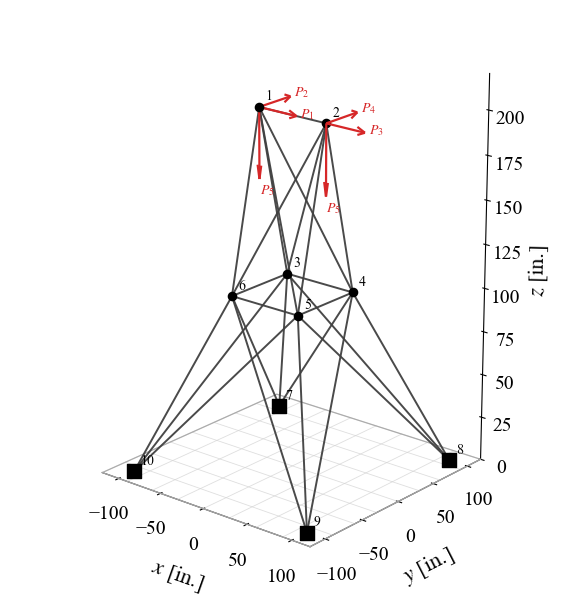}
\caption{The 25-bar transmission-tower-type space truss. The
structure contains 10 nodes and 25 members, with nodes 7--10
fully fixed. The uncertain loads $P_1$--$P_4$ act in the $x$-
and $y$-directions at nodes 1 and 2. The common vertical load
$P_5$ acts in the negative $z$-direction at both nodes.}
  \label{fig:truss25-schematic}
\end{figure}

The third benchmark is a 30-dimensional
compliance problem for a
25-bar transmission-tower-type space truss.
Figure~\ref{fig:truss25-schematic} illustrates the structural
geometry, supports, and loading. The structure contains 10 nodes
and 25 members, with nodes 7--10 fully fixed. This benchmark
tests SDD and GE-SDD in a high-dimensional, limited-sample regime
where the retained basis contains substantially more
coefficients than the number of function-value observations.

The input random vector is
$  \mathbf X
  =
  (A_1,\ldots,A_{25},
   P_1,P_2,P_3,P_4,P_5)^\top
  \in\mathbb R^{30},
  \label{eq:truss25-input}$
where $A_i$, $i=1,\ldots,25$, denotes the
cross-sectional area of the $i$-th member. The variables $P_1$ and $P_2$ define the $x$- and $y$-direction
loads at node 1, while $P_3$ and $P_4$ define the corresponding
loads at node 2. The variable $P_5$ defines the common load in
the negative $z$-direction at both nodes.
All input variables are mutually independent and
uniformly distributed.
Table~\ref{tab:truss25-inputs}
lists their bounds.

\begin{table}
\caption{Uncertain input variables and their probability
distributions for the 30-dimensional 25-bar space truss
compliance problem. All 30 input variables are mutually
independent and uniformly distributed as
$X_k\sim U(a_k,b_k)$.}
\label{tab:truss25-inputs}
\begin{tabular*}{\textwidth}{@{\extracolsep{\fill}}lllll@{}}
\toprule
Variable & Description & Unit & $a_k$ & $b_k$ \\
\midrule
$A_1,\ldots,A_{25}$
    & member cross-sectional areas
    & in$^2$
    & 0.5
    & 1.5 \\
$P_1$
    & $x$-direction load at node 1
    & lbf
    & $-2.0\times10^{3}$
    & $2.0\times10^{3}$ \\
$P_2$
    & $y$-direction load at node 1
    & lbf
    & $8.0\times10^{3}$
    & $1.2\times10^{4}$ \\
$P_3$
    & $x$-direction load at node 2
    & lbf
    & $-2.0\times10^{3}$
    & $2.0\times10^{3}$ \\
$P_4$
    & $y$-direction load at node 2
    & lbf
    & $8.0\times10^{3}$
    & $1.2\times10^{4}$ \\
$P_5$
    & vertical load at nodes 1 and 2
    & lbf
    & $-6.0\times10^{3}$
    & $-4.0\times10^{3}$ \\
\bottomrule
\end{tabular*}
\end{table}

The quantity of interest is the structural compliance
\begin{equation}
  Y(\mathbf X)
  =
  C(\mathbf X)
  =
  \mathbf F(\mathbf X)^\top
  \mathbf u(\mathbf X)
  =
  \mathbf F(\mathbf X)^\top
  \mathbf K(\mathbf X)^{-1}
  \mathbf F(\mathbf X),
  \label{eq:truss25-compliance}
\end{equation}
where $\mathbf{K}(\mathbf{X})\in\mathbb{R}^{18\times18}$,
$\mathbf{F}(\mathbf{X})\in\mathbb{R}^{18}$, and
$\mathbf{u}(\mathbf{X})\in\mathbb{R}^{18}$ denote the assembled global
stiffness matrix, nodal load vector, and nodal displacement vector,
respectively, defined on the $18$ free degrees of freedom that remain
after eliminating the constrained degrees of freedom at the fixed
nodes 7--10 (each of the 10 nodes carries three translational degrees
of freedom, and the 12 degrees of freedom at nodes 7--10 are removed).
The stiffness matrix depends on the member areas through the linear
assembly $\mathbf{K}(\mathbf{X})=\sum_{j=1}^{25}A_j\mathbf{K}_j$, where
$\mathbf{K}_j$ is the unit-area element stiffness contribution of
member~$j$. The load vector $\mathbf{F}(\mathbf{X})$ collects the
uncertain nodal loads $P_1,\ldots,P_5$: $P_1$ and $P_2$ act in the
$x$- and $y$-directions at node~1, $P_3$ and $P_4$ act in the $x$- and
$y$-directions at node~2, and $P_5$ acts in the negative
$z$-direction at both nodes~1 and~2. The Young's modulus is fixed at
$E=1.0\times10^7\,\mathrm{psi}$.

\paragraph{Analytical compliance sensitivities.}
The self-adjoint compliance formulation gives
\begin{align}
  \frac{\partial C}{\partial A_j}
  &=
  -\mathbf u^\top
  \frac{\partial\mathbf K}{\partial A_j}
  \mathbf u,
  \qquad
  j=1,\ldots,25,
  \label{eq:truss25-grad-area}
  \\
  \frac{\partial C}{\partial P_k}
  &=
  2
  \left(
    \frac{\partial\mathbf F}{\partial P_k}
  \right)^\top
  \mathbf u,
  \qquad
  k=1,\ldots,5.
  \label{eq:truss25-grad-load}
\end{align}
These expressions provide all 30 gradients without additional
state-equation solutions beyond the displacement analysis used
to evaluate the compliance.

\paragraph{Spline discretization.}
SDD and GE-SDD use quadratic B-splines with $p_k=2$ and maximum
interaction order $S=2$. The area dimensions use a finer
univariate basis than the load dimensions:
\begin{equation}
  n_k
  =
  \begin{cases}
    5, & k=1,\ldots,25,\\
    3, & k=26,\ldots,30.
  \end{cases}
  \label{eq:truss25-basis-resolution}
\end{equation}

Accordingly, the 25 area dimensions each contribute four
nonconstant univariate basis functions, whereas the five
load dimensions each contribute two.
The resulting number of retained basis functions is
\begin{align}
  n_\psi
  &=
  1
  +25(5-1)+5(3-1)
  +{25\choose2}(5-1)^2
  +25\cdot5(5-1)(3-1)
  +{5\choose2}(3-1)^2
  =
  5951
  \label{eq:truss25-basis-count}
\end{align}
functions. 

Standard SDD provides $M$ function-value equations
for estimating the 5951 coefficients.
GE-SDD provides one function-value equation and
30 derivative equations per training sample, producing
$31M$ stacked equations.
At $M=100$, standard SDD therefore uses 100 equations,
while GE-SDD uses 3100 function-gradient equations at
the same 100 input locations.
Although both systems remain underdetermined, the gradients
provide substantially more information for coefficient
estimation.

We test
$M\in\{100,150,200,300,500,1000\}
\label{eq:truss-training-sizes}$
and repeat each case over 20 independently generated training
sets.
The reference statistics and pointwise prediction errors
are evaluated using the same MCS sample of size $10^6$.
SDD and GE-SDD
estimate the mean and standard deviation directly from their
coefficients. Kriging estimates the same quantities from
predictions on the MCS sample.

\begin{figure}
  \centering
  \includegraphics[width=.99\textwidth]
  {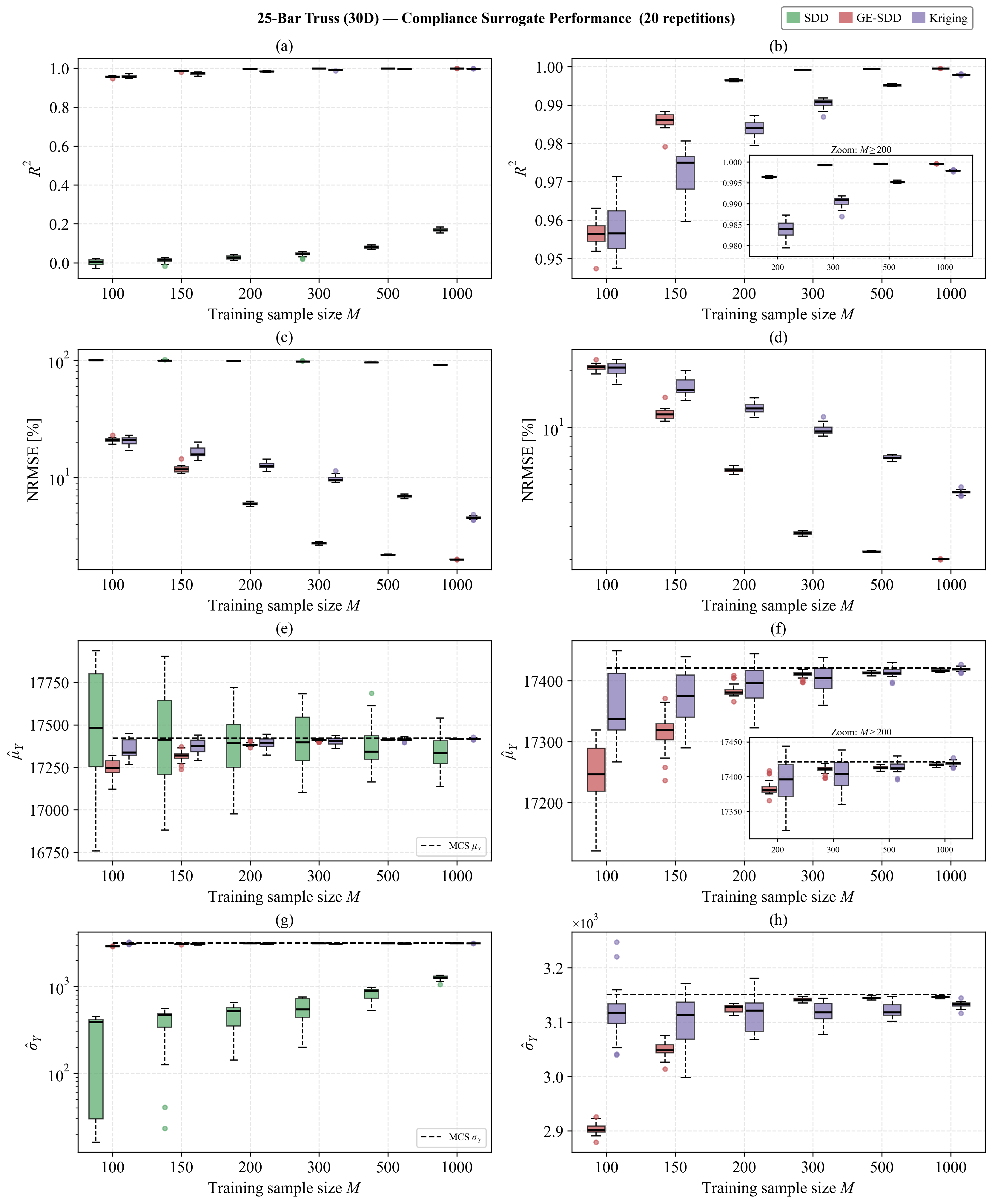}
    \caption{Prediction accuracy and moment estimates for the
    30-dimensional 25-bar space truss compliance problem over
    20 independent training-data realizations at
    $M\in\{100,150,200,300,500,1000\}$.
    The left column compares standard SDD, GE-SDD, and Kriging,
    whereas the right column provides an enlarged comparison
    between GE-SDD and Kriging.
    Panels~(a,~b) show the $R^2$ score, (c,~d) the NRMSE,
    (e,~f) the predicted mean $\hat{\mu}_Y$, and
    (g,~h) the predicted standard deviation $\hat{\sigma}_Y$.
    Dashed horizontal lines denote the MCS reference values.}
  \label{fig:truss25-boxplot}
\end{figure}

We compare GE-SDD with standard SDD and Kriging. We omit
GE-Kriging because the 30-dimensional problem would require
training on $31M$ correlated function and gradient observations. Repeated hyperparameter optimization at this scale is
computationally prohibitive for the present training-size sweep. 
Figure~\ref{fig:truss25-boxplot} summarizes the variability of
the prediction and moment-estimation metrics across replications.
Table~\ref{tab:truss25-results} reports their arithmetic means.
 
The left column of Figure~\ref{fig:truss25-boxplot}
shows that standard SDD remains inaccurate in pointwise
prediction across
all training sizes. 
With only $M$ function-value equations available for the
$5951$ retained coefficients, grouped cross-validation
 selects strong regularization. Because the ridge penalty excludes the constant
coefficient, SDD estimates the mean with average relative errors
between $0.63\%$ and $1.59\%$. However, regularization shrinks
the nonconstant coefficients and drives the surrogate toward a
constant-mean model. Consequently,  the average NRMSE of standard SDD decreases
only from $99.90\%$ at $M=100$ to $91.19\%$ at $M=1000$, while
its average $R^2$ increases only from $0.002$ to
$0.168$. SDD also underestimates the standard deviation, with
average relative errors decreasing from $91.1\%$ to $60.2\%$.
 
The accuracy comparison between GE-SDD and Kriging depends on
the training size. At $M=100$, the methods have similar
pointwise accuracy. The average NRMSE is $20.90\%$ for GE-SDD
and $20.58\%$ for Kriging, while the average $R^2$ is $0.9563$
and $0.9574$, respectively. GE-SDD produces the lower NRMSE in
8 of the 20 replications.
For every tested size from $M=150$ onward, GE-SDD achieves the lower NRMSE in
all 20 replications.  The average NRMSE is
$11.83\%$ for GE-SDD and $16.61\%$ for Kriging at $M=150$.
At $M=200$, the corresponding values are $5.95\%$ and
$12.69\%$. The difference increases for $M\ge300$. GE-SDD
decreases from $2.76\%$ at $M=300$ to $2.01\%$ at $M=1000$,
while Kriging decreases from $9.76\%$ to $4.54\%$. At $M=1000$,
the average $R^2$ reaches $0.9996$ for GE-SDD and $0.9979$ for
Kriging.

\clearpage

\begin{table*}
\caption{Mean prediction accuracy and moment-estimation
results for the 30-dimensional 25-bar space truss
compliance problem over 20 independent training-data
realizations.
Mean and standard-deviation errors are the averages over
the replications of the absolute relative errors with
respect to the MCS reference values
$\mu_{Y,\mathrm{MCS}}=1.7421\times10^{4}$ and
$\sigma_{Y,\mathrm{MCS}}=3.1509\times10^{3}$
($M_{\mathrm{MCS}}=10^6$).}
\label{tab:truss25-results}
\footnotesize
\begin{tabular*}{\textwidth}
{@{\extracolsep{\fill}}llrrrrrr@{}}
\toprule
$M$ & Method
    & NRMSE (\%, avg)
    & $R^2$ (avg)
    & $\hat{\mu}_Y$ (avg)
    & $\hat{\mu}_Y$ err.\ (\%, avg)
    & $\hat{\sigma}_Y$ (avg)
    & $\hat{\sigma}_Y$ err.\ (\%, avg) \\
\midrule
\multirow{3}{*}{100}
    & GE-SDD
    & 20.90
    & 0.9563
    & $1.725\times10^4$
    & 0.982
    & $2.904\times10^3$
    & 7.833 \\
    & SDD
    & 99.90
    & 0.002
    & $1.749\times10^4$
    & 1.589
    & $2.817\times10^2$
    & 91.061 \\
    & Kriging
    & 20.58
    & 0.9574
    & $1.736\times10^4$
    & 0.400
    & $3.118\times10^3$
    & 1.608 \\
\midrule
\multirow{3}{*}{150}
    & GE-SDD
    & 11.83
    & 0.9859
    & $1.732\times10^4$
    & 0.607
    & $3.050\times10^3$
    & 3.193 \\
    & SDD
    & 99.38
    & 0.012
    & $1.743\times10^4$
    & 1.410
    & $3.887\times10^2$
    & 87.664 \\
    & Kriging
    & 16.61
    & 0.9721
    & $1.737\times10^4$
    & 0.287
    & $3.103\times10^3$
    & 1.592 \\
\midrule
\multirow{3}{*}{200}
    & GE-SDD
    & 5.95
    & 0.9965
    & $1.738\times10^4$
    & 0.213
    & $3.125\times10^3$
    & 0.817 \\
    & SDD
    & 98.69
    & 0.026
    & $1.736\times10^4$
    & 0.946
    & $4.488\times10^2$
    & 85.755 \\
    & Kriging
    & 12.69
    & 0.9838
    & $1.739\times10^4$
    & 0.186
    & $3.114\times10^3$
    & 1.277 \\
\midrule
\multirow{3}{*}{300}
    & GE-SDD
    & 2.76
    & 0.9992
    & $1.741\times10^4$
    & 0.064
    & $3.141\times10^3$
    & 0.306 \\
    & SDD
    & 97.80
    & 0.043
    & $1.740\times10^4$
    & 0.784
    & $5.496\times10^2$
    & 82.557 \\
    & Kriging
    & 9.76
    & 0.9904
    & $1.741\times10^4$
    & 0.124
    & $3.118\times10^3$
    & 1.033 \\
\midrule
\multirow{3}{*}{500}
    & GE-SDD
    & 2.20
    & 0.9995
    & $1.741\times10^4$
    & 0.047
    & $3.144\times10^3$
    & 0.209 \\
    & SDD
    & 95.84
    & 0.081
    & $1.737\times10^4$
    & 0.693
    & $8.309\times10^2$
    & 73.631 \\
    & Kriging
    & 6.92
    & 0.9952
    & $1.741\times10^4$
    & 0.054
    & $3.122\times10^3$
    & 0.922 \\
\midrule
\multirow{3}{*}{1000}
    & GE-SDD
    & 2.01
    & 0.9996
    & $1.742\times10^4$
    & 0.025
    & $3.146\times10^3$
    & 0.155 \\
    & SDD
    & 91.19
    & 0.168
    & $1.733\times10^4$
    & 0.627
    & $1.255\times10^3$
    & 60.179 \\
    & Kriging
    & 4.54
    & 0.9979
    & $1.742\times10^4$
    & 0.019
    & $3.132\times10^3$
    & 0.600 \\
\bottomrule
\end{tabular*}
\end{table*}

The right column of Figure~\ref{fig:truss25-boxplot}
shows that GE-SDD variability decreases rapidly as $M$ increases. The
across-replication standard deviation of the GE-SDD NRMSE
decreases from $0.80$ percentage points at $M=100$ to less than
$0.01$ percentage points at $M=1000$. Kriging has greater
variability at every tested training size.
 
The moment estimates show a similar training-size transition.
Both methods estimate the mean accurately throughout the sweep.
Kriging has lower average relative mean errors for
$M\le200$. At $M=100$, $150$, and $200$, the Kriging errors are
$0.400\%$, $0.287\%$, and $0.186\%$, compared with
$0.982\%$, $0.607\%$, and $0.213\%$ for GE-SDD. The methods
produce similar mean estimates for $M\ge300$, and both average
relative errors fall below $0.03\%$ at $M=1000$. 

Kriging also estimates the standard deviation more accurately at
$M=100$ and $150$. GE-SDD becomes more accurate for
$M\ge200$. At $M=1000$, the average relative
standard-deviation error is $0.155\%$ for GE-SDD and $0.600\%$
for Kriging. Thus, gradient observations provide no clear
advantage at the smallest training size, where $31M$ remains
well below the 5951 retained coefficients. From $M=150$ onward,
GE-SDD provides more accurate pointwise predictions, and from
approximately $M=200$ onward, it provides more accurate
standard-deviation estimates.

\subsection{Discussion}
\label{subsec:discussion}

GE-SDD uses $M$ function observations and $M\times N$
partial-derivative observations.
The number of distinct input  locations therefore remains $M$,
but the cost of obtaining the gradients depends on the
underlying model.
Example~1 provides explicit derivatives. 
In Example~2, direct sensitivity analysis requires
additional linear-system solutions for the three uncertain
parameters, but all solutions reuse the same dynamic-stiffness factorization.
In Example~3, the self-adjoint compliance formulas require no additional state-equation solutions. 
For more general applications, adjoint sensitivity
analysis or automatic differentiation can reduce the
marginal cost of gradient observations relative
to additional high-fidelity simulations.

The algebraic cost of GE-SDD depends primarily on the
retained basis size
$n_\psi=L_{S,\mathbf p,\boldsymbol\Xi}$ and the $M(1+N)$ stacked equations.
A dense Cholesky factorization of the
derivative-augmented spline monomial moment matrix has computational
complexity $\mathcal O(n_\psi^3)$.
Consequently, the cost increases with input dimension,
interaction order, and basis resolution.
Large-scale applications can reduce them through sparse linear algebra, iterative solvers, dimension-adaptive basis selection, and further use of the tensor-product and ANOVA structures.

The accuracy gains reported in
Table~\ref{tab:example1-results},
Table~\ref{tab:2dof}, and
Table~\ref{tab:truss25-results} depend on the benchmark, training budget, basis resolution, and gradient availability.
Example~3 shows that gradient augmentation does not guarantee greater accuracy at every training size.
GE-SDD and Kriging have similar pointwise accuracy at $M=100$. GE-SDD becomes more accurate at $M\ge150$, and the difference increases at $M\ge300$.

The three benchmarks show when gradient augmentation is most beneficial.
In Example~2 with $N=3$ and $M=919$, Kriging remains the most accurate method because it re-estimates the anisotropic length scales at each frequency. GE-SDD is most accurate for the nonsmooth response in Example~1 and for the high-dimensional problem in Example~3 at moderate and large training sizes. In these settings, case-specific maximum-likelihood adaptation is less effective or computationally impractical. The observed crossover training sizes are not universal; they depend on the retained basis size, regularization, response structure, and the accuracy and cost of the gradients. Further studies using large-scale finite-element and computational fluid dynamics models with solver-based sensitivities are needed to assess the broader accuracy and computational efficiency of GE-SDD.

\section{Conclusions}
\label{sec:conclusions}

This paper proposed a gradient-enhanced spline
dimensional decomposition (GE-SDD) method for uncertainty
quantification (UQ) with limited training samples, including
high-dimensional regimes. GE-SDD
combines function values and partial derivatives through
a probability-weighted Sobolev transformation, balances
the function and derivative blocks, and estimates the coefficients
using ridge regression with training-sample-grouped $K$-fold cross-validation. A final coefficient-estimation step transforms the solution to the original $L^2$-orthonormal SDD coordinates, preserving closed-form mean and variance estimates.

Numerical results from three benchmark problems showed that
GE-SDD is more accurate than standard SDD and GE-Kriging, while its comparison with function-only Kriging depends on the problem.
For the two-dimensional continuous exponential function, GE-SDD
achieved a median NRMSE of $0.269\%$, compared with $8.086\%$
for Kriging, the most accurate baseline in that benchmark.
For the 2-DOF dynamical system, GE-SDD reduced the MAE of the
predicted frequency-dependent standard deviation to
$2.435\times10^{-5}$, compared with $4.978\times10^{-5}$ for
standard SDD and $4.651\times10^{-5}$ for GE-Kriging. GE-SDD was more accurate than both in every replication with the smallest
across-replication variability; the frequency-wise
maximum-likelihood-tuned Kriging baseline attained the lowest
MAE of $1.151\times10^{-5}$ in this low-dimensional benchmark.
For the 30-dimensional truss, GE-SDD and Kriging had similar pointwise accuracy at $M=100$. GE-SDD produced lower NRMSE and higher $R^2$ for $M\ge150$, with a larger difference at moderate and large training sizes.

The truss results also showed that the relative accuracy depended on the target UQ quantity. Both methods estimated the mean accurately across all training sizes, although Kriging was slightly more accurate at smaller sizes. GE-SDD provided more accurate standard-deviation estimates for $M\ge200$. Conversely, Example~2 showed that anisotropic Kriging could remain more accurate near sharp resonances when the input dimension was small and frequency-specific hyperparameter optimization was computationally feasible. The benefits of gradient augmentation therefore depended on the input dimension, training size, basis resolution, response structure, and target statistic.

Future work will evaluate GE-SDD on large-scale finite-element and computational fluid dynamics models using adjoint sensitivities or automatic differentiation. Further developments will integrate GE-SDD with first- and second-order reliability methods, introduce gradient-guided adaptive knot placement, and develop dimension-adaptive basis selection for responses with localized rapid variations such as resonance zones.



\section*{Funding}
This work was supported by the National Research Foundation of Korea (NRF) grant funded by the Korea government (MSIT) (RS-2025-00560781).

\section*{Declaration of generative AI and AI-assisted
technologies in the writing process}

During the preparation of this work, the authors used 
OpenAI's ChatGPT and Anthropic's Claude for English-language
editing, grammar checking, wording refinement, and stylistic
polishing. After using these tools, the authors reviewed and
edited the content as needed and take full responsibility 
for the content of the published article.

\appendix
\section{Construction of measure-consistent orthonormal
B-spline basis functions}
\label{app:basis-construction}
This appendix presents the detailed construction of the
measure-consistent orthonormal B-spline basis functions
used in the SDD representation.
\begin{itemize} 
\item[(1)] Specify the univariate spline spaces.
For each input coordinate $X_k$, let $p_k$, $\boldsymbol{\xi}_k$, and $n_k$ denote the B-spline degree, open knot vector, and number of basis functions, respectively. Collect these parameters as
\begin{equation}
\mathbf{p}:=(p_1,\ldots,p_N),
\qquad
\mathbf{n}:=(n_1,\ldots,n_N),
\qquad
\boldsymbol{\Xi}
:=
(\boldsymbol{\xi}_1,\ldots,\boldsymbol{\xi}_N).
\label{eq:spline-parameter-collections}
\end{equation}
For \(u\subseteq\{1,\ldots,N\}\), let \(\mathbf{p}_u\), \(\mathbf{n}_u\), and \(\boldsymbol{\Xi}_u\) denote the corresponding subcollections.
\item[(2)] Construct the raw univariate bases.
For the \(k\)-th input coordinate, define
\begin{equation}
  \mathbf{P}_k(X_k)
  :=
  \bigl(
    1,\,
    B^k_{2,p_k,\boldsymbol{\xi}_k}(X_k),\,
    \ldots,\,
    B^k_{n_k,p_k,\boldsymbol{\xi}_k}(X_k)
  \bigr)^\top,
  \label{eq:Pvec}
\end{equation}
$B^k_{i,p_k,\boldsymbol{\xi}_k}$ is the $i$-th B-spline of degree $p_k$, generated by the Cox--de~Boor recursion~\cite{Cox1972,deBoor1972}. Replacing the first B-spline with $1$ preserves the constant direction during orthonormalization and produces nonconstant basis functions with zero mean under $f_{X_k}(\mathbf{X}_k)$.
\item[(3)] Orthonormalize the univariate bases. Define the Gram matrix 
\begin{equation}
  \mathbf{G}_k
  :=
  \mathbb{E}\!\left[
    \mathbf{P}_k(X_k)\mathbf{P}_k^\top(X_k)
  \right]
  =
  \int_{a_k}^{b_k}
  \mathbf{P}_k(x)\mathbf{P}_k^\top(x)
  f_{X_k}(x)\,\mathrm{d}x
  \in\mathbb{R}^{n_k\times n_k}.
  \label{eq:gram}
\end{equation}
We evaluate this integral over each nonzero knot span using 24-point Gauss--Legendre quadrature. We use the same rule for the basis-derivative products required by GE-SDD.
Because $\mathbf{G}_k$ is symmetric positive definite, its Cholesky factorization is
$
\mathbf{G}_k=
\mathbf{Q}_k\mathbf{Q}_k^\top,
\label{eq:gram-cholesky-univariate}
$
where $\mathbf{Q}_k$ is nonsingular and lower triangular. The measure-consistent orthonormal basis is
\begin{equation}
\boldsymbol{\psi}_k(X_k)
:=
\mathbf{Q}_k^{-1}\mathbf{P}_k(X_k)
=
\left(
\psi^k_{1,p_k,\boldsymbol{\xi}_k}(X_k),
\ldots,
\psi^k_{n_k,p_k,\boldsymbol{\xi}_k}(X_k)
\right)^\top.
\label{eq:ortho-univariate}
\end{equation}
This basis satisfies
$
\mathbb{E}[\boldsymbol{\psi}_k(X_k)]
=(1,0,\ldots,0)^\top,
~\text{and}~
\mathbb{E}\left[
\boldsymbol{\psi}_k(X_k)
\boldsymbol{\psi}_k^\top(X_k)
\right]=
\mathbf{I}_{n_k}.
\label{eq:ortho-char}
$
Thus, $\psi^k_{1,p_k,\boldsymbol{\xi}_k}\equiv1$, whereas
$
\mathbb{E}\left[
\psi^k_{i,p_k,\boldsymbol{\xi}_k}(X_k)
\right]=
0,
~ i=2,\ldots,n_k.
\label{eq:univariate-zero-mean}
$
\item[(4)] Form the tensor-product bases.
For a non-empty interaction subset
\begin{equation}
  u=\{k_1,\ldots,k_{|u|}\}
  \subseteq\{1,\ldots,N\},
\end{equation}
define the full multi-index set
\begin{equation}
  \mathcal{I}_{u,\mathbf{n}_u}
  :=
  \left\{
    \mathbf{i}_u= (i_{k_1}, \dots, i_{k_{|u|}}) :
    1\le i_k\le n_k,\; k\in u
  \right\}.
  \label{eq:full-index-set}
\end{equation}
For
$\mathbf{i}_u\in\mathcal{I}_{u,\mathbf{n}_u}$,
the corresponding tensor-product basis function is
\begin{equation}
  \Psi^u_{\mathbf{i}_u,\mathbf{p}_u,\boldsymbol{\Xi}_u}
  (\mathbf{x}_u)
  :=
  \prod_{k\in u}
  \psi^k_{i_k,p_k,\boldsymbol{\xi}_k}(x_k).
  \label{eq:tp-basis-appendix}
\end{equation}
\end{itemize} 
The full index set in \eqref{eq:full-index-set} includes the constant index $i_k=1$. If one or more active coordinates take this constant index, the corresponding tensor-product basis function effectively reduces to a lower-order interaction. Therefore, to prevent the degree of interaction from falling below $|u|$, the retained SDD basis excludes the constant index in every active coordinate.
\bibliographystyle{elsarticle-num}
\bibliography{citation}

\end{document}